\documentclass[onefignum,onetabnum,amsmath,amssymb,twocolumn]{article}

\usepackage{graphicx}
\usepackage{dcolumn}
\usepackage{bm}

\usepackage{amssymb}
\usepackage{amsfonts}
\usepackage{amsmath}
\usepackage{euscript}
\usepackage{algorithm}
\usepackage{algpseudocode}
\usepackage{mathrsfs}
\usepackage{hyperref}
\usepackage{color}
\usepackage[margin={1.2cm,2cm}]{geometry}

\newtheorem{theorem}{Theorem}
\newtheorem{corollary}{Corollary}
\newtheorem{remark}{Remark}

\newcommand{\bo}[1]{\mathbf{#1}}
\newcommand{\bb}[1]{\mathbb{#1}}
\newcommand{\cc}[1]{\mathcal{#1}}
\newcommand{\kk}[1]{\mathfrak{#1}}
\newcommand{\scr}[1]{\mathscr{#1}}

\newcommand{\nc}{\newcommand}
\nc{\cblue}[1]{\color{blue}{#1}\color{black}{}}
\nc{\cred}[1]{\color{red}{#1}\color{black}{}}
\nc{\cmblue}[1]{\color{MidnightBlue}{#1}\color{black}{}}
\nc{\cgrey}[1]{\color{gray}{#1}\color{black}{}}
\nc{\cmagenta}[1]{\color{magenta}{#1}\color{black}{}}
\nc{\corange}[1]{\color{orange}{#1}\color{black}{}}
\nc{\cviolet}[1]{\color{violet}{#1}\color{black}{}}

\begin{document}

\title{ANKH: A Generalized $\cc{O}(N)$ Interpolated Ewald  Strategy for Molecular Dynamics Simulations}

\author{Igor Chollet\thanks{LAGA, Université Sorbonne Paris Nord, UMR 7539, Villetaneuse, France and LCT, Sorbonne Université, UMR 7616, Paris, France. chollet@math.univ-paris13.fr}
\and Louis Lagardère\thanks{LCT, Sorbonne Université, UMR 7616, Paris, France. louis.lagardere@sorbonne-universite.fr}
\and Jean-Philip Piquemal\thanks{LCT, Sorbonne Université, UMR 7616, Paris, France. jean-philip.piquemal@sorbonne-universite.fr}
}
\date{\today}

\maketitle
\begin{abstract}
To evaluate electrostatics interactions, Molecular dynamics (MD) simulations rely on Particle Mesh Ewald (PME), an $\cc{O}(Nlog(N))$ algorithm that uses Fast Fourier Transforms (FFTs) or, alternatively, on $\cc{O}(N)$ Fast Multipole Methods (FMM) approaches. However, the FFTs low scalability remains a strong bottleneck for large-scale PME simulations on supercomputers. On the opposite, - FFT-free - FMM techniques are able to deal efficiently with such systems but they fail to reach PME performances for small- to medium-size systems, limiting their real-life applicability. We propose ANKH, a strategy grounded on interpolated Ewald summations and designed to remain efficient/scalable for any size of systems.
The method is generalized for distributed point multipoles and so for induced dipoles which makes it suitable for high performance simulations using new generation polarizable force fields towards exascale computing.
\end{abstract}

\section{Introduction}
\label{s::intro}
Performing Molecular Dynamics (MD) simulations requires to numerically solve Newton's equations of motion for a system of interacting particles. Modern simulations rely on Molecular Mechanics to evaluate the different physical interactions acting on an ensemble of particles. These approaches are very diverse and range from classical force fields (FFs)  \cite{case2005amber,brooks2009charmm,van2005gromacs} to more evolved polarizable approach embodying many-body effects (PFFs) \cite{chapterpol2015,Meclr19,annurev-biophys}. In any case, one needs to compute electrostatic interactions that are associated to the Coulomb energy and are an essential contribution to the systems total potential energy. To efficiently compute these quantities, MD softwares mainly exploit grid based methods such as the Particle Mesh Ewald (PME) approach, which is an $\cc{O}(Nlog(N))$ algorithm that uses Fast Fourier Transforms (FFTs). While PME is extremely efficient for small to medium system sizes, for very large ensembles of particles, the FFTs limited scalability is a major bottleneck for large-scale simulations on supercomputers. Historically, Fast Multipole Methods (FMM) have been considered as good $\cc{O}(N)$ candidates to overcome such limitations.  The FMM strategy is - FFT-free and capable of efficiently dealing with such very large systems. Still, PME performance is yet to be reached for small- to medium-size systems, limiting the real-life applicability of FMMs.
In this paper, we introduce ANKH, a new general strategy for the fast evaluation of the electrostatic energy using periodic boundary conditions and a density of charge due to point multipoles:
\begin{equation}
\label{eq_qmutheta_coulomb}
    \cc{E}:=\frac{1}{2}\sum_{\bo{t}\in 2r_B\bb{Z}^3}\sum_{\bo{x}\in X}\sum_{\bo{y}\in Y}\kk{D}_{\bo{x}}\kk{D}_{\bo{y}}\left(\frac{1}{|\bo{x}-\bo{y}+\bo{t}|}\right),
\end{equation}
where $X,Y\subset B\subset \bb{R}^3$ are two point clouds (of atoms) that can be equal, $\kk{D}_{u} := q_{\bo{u}} + \mu_{\bo{u}}\cdot \nabla + \Theta_{\bo{u}}:\nabla^2$, $q_\bo{u}$ is the charge of the atom $\bo{u}$, $\mu_\bo{u}$ its dipole moment, $\Theta_\bo{u}$ its quadrupole moment and $\nabla$, $\nabla^2$ respectively denote the gradient and Hessian operators, $B$ is the simulation box of radius $r_B$. Each point cloud $X$ and $Y$ is composed of $N$ atoms, named \textit{particles} to fit with the usual notation used in hierarchical methods \cite{GREENGARD1987325}. Notice that we restricted ourselves to charges, dipoles and quadrupoles for the sake of clarity, but that all the theory presented in this paper can be straightforwardly extended to any multipole order and therefore to induced dipoles, providing solutions for polarizable force fields.

Naive computation of the energy through direct implementation of Eq. \eqref{eq_qmutheta_coulomb} faces numerical convergence issues, requiring reformulation. The literature widely exploits Ewald summation techniques \cite{pmeref,spmeref,toukmaji2000efficient,sagui2004towards,Stamm_2018} to provide numerically convergent expressions to compute Eq. \eqref{eq_qmutheta_coulomb}. 
The previously mentioned limitations of 3D-FFT based techniques motivate the development of alternatives \cite{doi:10.1063/5.0040966,dardenbrooks,doi:10.1021/ct5009075}, also based on Ewald summation, such as fast summation schemes for Eq. \eqref{eq_qmutheta_coulomb} directly \cite{OHNO20142575,LAMBERT1996274,doi:10.1021/acs.jctc.0c00744} and applied to molecular dynamics, mainly exploiting Fast Multipole Method \cite{GREENGARD1987325} (FMM) and cutoff approaches \cite{ewaldmultipole}.
The latter has many advantages since it provides linear complexity and high scalability, despite the error that may occur when applying cutoff on images of $B$. Convergent alternatives \cite{Schmidt1991ImplementingTF,challacombe} also extend FMM for Coulomb potential to periodic case when $\kk{D}_\bo{x}$ and $\kk{D}_{\bo{y}}$ are restricted to charges. Other kernels also benefit from periodic extensions \cite{407723}.

Our ANKH approach aims at providing a theoretical framework well-suited for linear-complexity and scalable FMM-based energy computations, built on Ewald summation. The methodology presented here introduces various novelties, declined in two variants, both exploiting different ways of solving $N$-body problems appearing in Ewald summation. Among them, we introduce:
\begin{itemize}
    \item a new interpolated Ewald summation, leading to numerical schemes to handle differential operators of multipolar (polarizable) molecular dynamics, accelerated through two different numerical methods,
    \item alternative techniques to account for periodicity,
    \item explicit formulae to handle the mutual interactions.
\end{itemize}

Due to the positioning of our work, this article deals with mathematical and computer science topics in the precise context of molecular dynamics. For the sake of clarity, we provide both the theory and the algorithms that should allow to minimize the reader's effort to implement our method. However, efficient implementations rely on various optimizations and some details are beyond the scope of this paper. This is why we also provide our code, used to generate the results in the following.

The paper is organized as follows. We first briefly review in Sect. \ref{s::ibfmm} the interpolation-based Fast Multipole Method (referred to as IBFMM in the remainder), then we introduce in Sect. \ref{s::ankhfmm} ANKH-FMM which is our new method for the fast Ewald Summation. In Sect. \ref{s:ankhfft} we introduce an entirely new approach named ANKH-FFT that overcomes some limitations of ANKH-FMM regarding modern High Perfoamnce Computing (HPC) architectures (at the cost of a slightly higher complexity than ANKH-FMM), then we present in Sect. \ref{s::numerical} numerical results to emphasize on the performance of our method as well as a comparison with a production code for molecular dynamics (namely Tinker-HP  \cite{lagardere2018tinker}). 

\section{Interpolation-based FMM}
\label{s::ibfmm}
This section is dedicated to the presentation of interpolation-based FMM. We first recall in Sect. \ref{ss::generalities} the reasons of the efficiency of the original 3D FMM (for Coulomb potential) as well as the algorithm in the case of quasi-uniform particle distributions. Then, in Sect. \ref{ss::ibfmmkis}, we shortly describe one of the main kernel-independent FMM, namely the interpolation-based FMM, providing the corresponding operator formulas.
\subsection{Generalities of FMMs}
\label{ss::generalities}
Fast Multipole Methods \cite{GREENGARD1987325,10.1007/BFb0089775,FONG20098712} (FMM) is a family of fast divide and conquer strategies to compute $N$-body problems in linear or linearithmic complexity. Important efforts have been made in the literature to incorporate 3D FMM (i.e. for the Coulomb potential) in molecular dynamics simulation codes \cite{shimada,doi:10.1021/acs.jctc.0c00744,LAMBERT1996274}. Such a FMM scheme is based on a separable spherical harmonic expansion of the Coulomb potential $\frac{1}{|\bo{x}-\bo{y}|}$ writing as:
\begin{equation}
\label{coulombfmm}
    \sum_{l=0}^{+\infty}\frac{4\pi r^l_<}{(2l+1)r^{l+1}_>}\sum_{m=-l}^lY^m_l(\theta_\bo{x},\phi_\bo{x})Y^{-m}_l(\theta_\bo{y},\phi_\bo{y}),
\end{equation}
where $\bo{x} = (r_\bo{x},\theta_\bo{x},\phi_\bo{x})$, $\bo{y} = (r_\bo{y},\theta_\bo{y},\phi_\bo{y})$ in spherical coordinates, $r_<:= min\{r_\bo{x},r_\bo{y}\}$, $r_>:= max\{r_\bo{x},r_\bo{y}\}$ and $Y^m_l$ refers to the spherical harmonic of order $(l,m)$. Restricting to subset of particles $\bo{x}\in t\subset\bb{R}^3,\bo{y}\in s\subset\bb{R}^3$ such that $\frac{radius(t)+radius(s)}{distance(t,s)}\leq \nu$, $\nu>0$ sufficiently small, one may approximate Eq. \eqref{coulombfmm} by truncating this series to a finite order $P$. This allows to derive a fast summation scheme for any $\bo{x}\in X$, for problems consisting in computing $\phi$ such that 
\begin{equation}
\label{nbodypb}
    \phi(\bo{x}) := \sum_{\bo{y}\in Y}\frac{q(\bo{y})}{|\bo{x}-\bo{y}|}
\end{equation}
where $q(\bo{y})$ is a charge associated to the particle $\bo{y}$. Such a problem in Eq. \eqref{nbodypb} is referred to as a $N$-body problem and such $t$ and $s$ are named \textit{well-separated target ($t$) and source ($s$) cells}, where a \textit{cell} denotes in this paper a cubical subset of $\bb{R}^3$. These cells are obtained in practice through a recursive splitting of the computational box $B$ into smaller boxes. Assuming that $B$ is a cube, one may divide it into $8$ other cubes and repeat this procedure until the induced cubes (the cells) are sufficiently small. A tree whose nodes are cells and such that the daughters of a given cell are the non-empty cubes coming from its decomposition is named an \textit{octree}. The root cell is considered to belong at octree level $0$ and any non-root cell belongs at level $E+1$, where $E$ denotes its mother level. In the case of perfect octrees (i.e. with maximal amount of daughter per non-leaf node, that is $8$ daughters), we consider that two cells at the same tree level with strictly positive distance (i.e. non-neighbor cells) are well-separated.

\begin{figure}
    \centering
    \includegraphics[width=\linewidth]{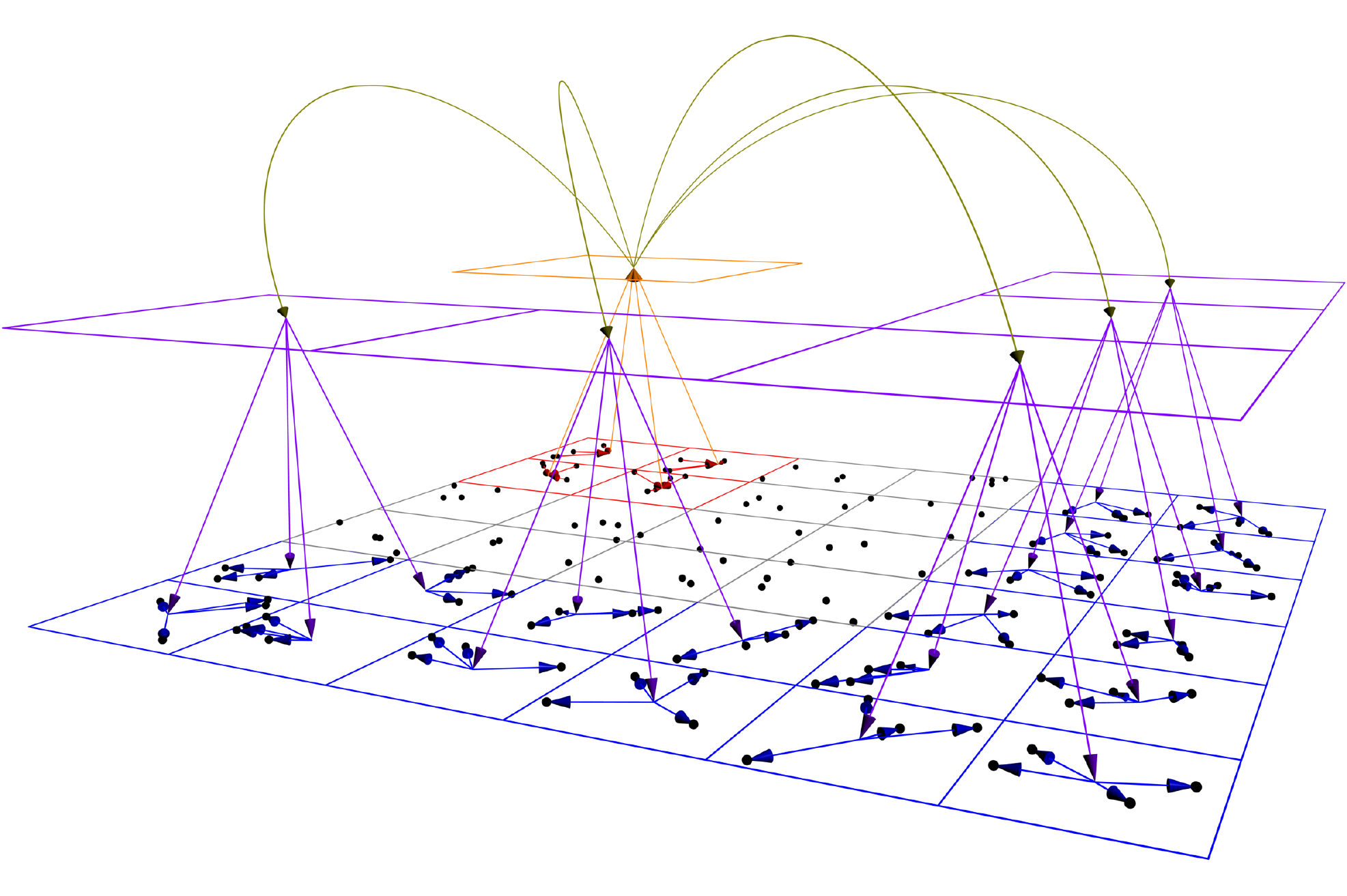}
    \caption{2D schematic representation of FMM operators. A source cell $s$ at level 2 (orange) interacts with well-separated target cells $t$'s at level 2 (purple) through M2L operators (green arrows). At level 3, daughters (red) of $s$ aggregate their multipole expansions through M2M operators (orange arrows) to form multipole expansion in $s$. Each $t$ generates 4 L2L operators to add influence of $s$ to its daughters through L2L operators (purple arrows). Red and blue arrows in source and target cells respectively represent the action of the P2M and L2P operators respectively in the corresponding cell. Grey cells at level 3 correspond to the daughter of non well-separated cells from $s$ at level 2.\label{fig:fmmoperators}}
\end{figure}

The sum over $\bo{y}\in Y_{|s}$ (the restriction of $Y$ to the particles in $s$) can be switched, for any $\bo{x}\in X_{|t}$ to the terms depending only on $\bo{y}$, due to the separability in the variables $\bo{x}$ and $\bo{y}$ of Eq. \eqref{coulombfmm}. Expression in  Eq. \eqref{coulombfmm} results \cite{10.1007/BFb0089775} in a three-step summation scheme expressing Eq. \eqref{nbodypb} in the form 
\begin{equation}
\label{eq:generalfmm}
    \sum_kU_k[t]^*(\bo{x})\sum_lC_{k,l}[t,s]\sum_{\bo{y}\in Y_{|s}}V_l[s](\bo{y})q(\bo{y})
\end{equation}
or equivalently under matrix form $U^*[t]C[t,s]V[s]q_{|s}$ where $q_{|s}$ refers to the vector of charges of particles in $s$. The result $\cc{M}_s$ of the product $V[s]q_{|s}$ is named the \textit{multipole expansion} of source cell $s$. The matrix $V[s]$ involved in this product is named Particle-To-Multipole operator (shortly denoted by P2M) since it maps information associated to particles (i.e. here the charges associated to atoms) to the multipole expansion of the cell in which they belong. This multipole expansion is transformed into a \textit{local} one $\cc{L}_t$ by means of product by the Multipole-To-Local operator (shortly denoted by M2L) $C[t,s]$. Similarly, the product by $U^*$ is referred to as the Local-To-Particle operator (L2P). Nevertheless, the efficiency of the FMM algorithm is also guaranteed by two additional operators based on the idea that multipole expansions of a given non leaf cell can be approximating through combination of its daughter's multipole expansions by means of a Multipole-To-Multipole operator (M2M) and local expansions in a given non root cell can be approximated from its mother local expansion thanks to a Local-To-Local operator (L2L). Finally, two leaf-cells that are not well-separated still need to interact, but without expansion. This is done by applying the Particle-To-Particle operator (P2P), i.e. nothing else than a direct computation involving only the particles of these two cells. The relation between all these operators according to the tree structure are depicted in Fig. \ref{fig:fmmoperators}.

In the case of quasi-uniform particle distributions (i.e. with approximately the same density of particles in all the simulation box), octrees can be considered as perfect. In this situation, the FMM algorithm is easy to provide (and usually referred to as the \textit{uniform FMM algorithm} in the literature), according to the definition of the different FMM operators. This is summarized in Alg. \ref{alg:uniffmm}, where we used the notation $V[s,s']$ to denote the M2M operator between the source cell $s$ and one of its daughter $s'$; $U^*[t',t]$ to denote the L2L operator between the target cell $t$ and one of its daughter $t'$; for any target cell $t$, $\Lambda(t)$ the interaction list of $t$, that is the set of source cells $s$ that are well-separated from $t$ and whose ancestors in the octree are not well-separated from any ancestor of $t$; $\scr{T}$ refers to the octree; $\cc{L}(\scr{T})$ refers to the set of leaves of $\scr{T}$ and $K(\bo{x},\bo{y}) = \frac{1}{|\bo{x}-\bo{y}|}$.
\begin{algorithm}[htbp]
\caption{Uniform FMM algorithm}
\label{alg:uniffmm}
\begin{algorithmic}
    \State // \textbf{\underline{Step 1:} compute multipole expansions}
    \For{each leaf cell $s\in \cc{L}(T)$}
        \State // Compute multipole expansion $\cc{M}_s$ in $s$ using P2M
        \State $\cc{M}_s = V[s]q_{|s}$
    \EndFor
    \For{each non leaf level $e$ starting from the deepest one}
        \For{each cell $s$ at level $e$}
            \State // Aggregate multipole expansions through M2M
            \State $\cc{M}_s = \bo{0}$
            \For{each daughter $s'$ of $s$}
                \State $\cc{M}_s = \cc{M}_s + V[s,s']\cc{M}_{s'}$
            \EndFor
        \EndFor
    \EndFor
    \State // \textbf{\underline{Step 2:} Compute local expansions}
    \For{each target cell $t$}
        \State // Add far contribution using M2L
        \State $\cc{L}_t = \bo{0}$
        \For{each source cell $s\in \Lambda(t)$}
            \State $\cc{L}_t = \cc{L}_t+C[t,s]\cc{M}_s$
        \EndFor
    \EndFor
    \For{each level $e$ starting from the root}
        \State // Scatter local expansions using L2L
        \For{each target cell $t$ at level $e$}
            \For{each daughter $t'$ of $t$}
                \State $\cc{L}_t = \cc{L}_t + U^*[t',t]\cc{L}_t$
            \EndFor
        \EndFor
    \EndFor
    \State // \textbf{\underline{Step 3:} convert local expansions}
    \For{each leaf cell $t\in \cc{L}(\scr{T})$}
        \State // Apply L2P
        \State $\phi_t = U^*[t]\cc{L}_t$
    \EndFor
    \State // \textbf{\underline{Step 4:} add near contribution}
    \For{each leaf cell $t\in \cc{L}(\scr{T})$}
        \State // Apply all P2P involving $t$
        \For{each source cell $s$ adjacent to $t$}
            \For{$\bo{x}\in X_{|t}$}
                \State $\phi(\bo{x}) = \phi_t(\bo{x})$
                \For{$\bo{y}\in Y_{|s}$}
                    \State $\phi(\bo{x}) = \phi(\bo{x}) + K(\bo{x},\bo{y})q(\bo{y})$
                \EndFor
            \EndFor
        \EndFor
    \EndFor
\end{algorithmic}
\end{algorithm}

Alg. \ref{alg:uniffmm} results in $\cc{O}(N)$ complexity \cite{10.1007/BFb0089775} with prefactor depending on the user required precision. Notations of Alg. \ref{alg:uniffmm} directly provide a way of approximating the influence $\phi_{t',s'}$ of $s'$ a daughter of $s$ on $t'$ and daughter of $t$, provided that $t$ and $s$ are well-separated, under matrix form:
\begin{equation}
\label{eq:fmmchain}
    \phi_{t',s'} \approx \underbrace{U^*[t']}_{L2P}\underbrace{U^*[t',t]}_{L2L}\underbrace{C[t,s]}_{M2L}\underbrace{V[s,s']}_{M2M}\underbrace{V[s']}_{P2M}q_{|s'}.
\end{equation}
This last can be interpreted as low-rank matrix approximation when seeing Eq. \eqref{nbodypb} as a matrix-vector product with matrix entries corresponding to evaluation of the Coulomb potential on pairs of target and source cells. Notice that the interaction between equally located particles is not well-defined, and we consider in this article that they are just equal to $0$.

In terms of parallel implementation, 3D FMM as shown its ability to scale in distributed memory context. Combined with its linear complexity, this makes this method particularly attractive for MD simulations. However, few limitations still remain: the presentation we made only takes into accounts charges and considered non-periodic case. This last can be solved by means of periodic FMM \cite{challacombe}.

\subsection{Interpolation-based kernel-independent scheme}
\label{ss::ibfmmkis}
Actually, a similar scheme can be derived for any kernel $K(\bo{x},\bo{y})$ (not only the Coulomb one) that satisfies particular assumption, namely the \textit{asymptotically smooth} behavior \cite{bookborm}. There exist different approaches to obtain a fast separable formula in a \textit{kernel-independent} way \cite{YING2004591,FONG20098712}. Among them, we are especially interested into the interpolation-based FMM (IBFMM) because of its flexibility, its well documented error bounds as well as the particular form of the approximated expansion of $K$ it provides. Indeed, denoting by $S_k[t]$ (resp. $S_l[s]$) a multivariate Lagrange polynomial associated to the interpolation node $\bo{x}_k^t$ (resp. $\bo{y}_l^s$) in $t$ (resp. $s$), a Lagrange interpolation of $K$ in $t\times s$ writes:
\begin{equation}
\label{interpolationformula}
    K(\bo{x},\bo{y}) \approx \sum_{k}S_k[t](\bo{x})\sum_lK(\bo{x}_k^t,\bo{y}_l^s)S_l[s](\bo{y})
\end{equation}
for any $(\bo{x},\bo{y})\in t\times s$. Eq. \eqref{interpolationformula} clearly exhibits a separable expansion of $K$ in the $\bo{x}$ and $\bo{y}$ variables. Hence, the restriction of Eq. \eqref{nbodypb} to particles in $t\times s$ (i.e. to $X_{|t}$ and $Y_{|s}$) and replacing the Coulomb kernel by a generic one $K$) can be approximated by:
\begin{equation}
\underbrace{\sum_{k}S_k[t](\bo{x})}_{\text{L2P on }t}\underbrace{\sum_lK(\bo{x}_k^t,\bo{y}_l^s)}_{\text{M2L between $t$ and $s$}}\underbrace{\sum_{\bo{y}\in Y_{|s}}S_l[s](\bo{y})q(\bo{y})}_{\text{P2M on $s$}}
\end{equation}
by analogy with Eq. \eqref{eq:generalfmm}. This identifies most of the FMM operators. Moreover, interpreting multipole expansions generated this way as charges associated to interpolation nodes (now seen as particles), one can obtain a new $N$-body problem with the same kernel $K$. Hence, this procedure can be repeated all along the tree structure in order to approximate multipole expansions in the non-leaf cells (which is the analogous of the M2M operator). In the same way, evaluating interpolant at a non-leaf cell on interpolation nodes in its daughters allows to derive a analogous of the L2L operator. To be more precise, we assume that the Lagrange interpolation applied on  cell $c$ uses the interpolation grid $\Xi_c=\{\bo{y}^s_0,...,\bo{y}^s_{P-1}\}$ with $P$ interpolation nodes over $c$. In practice, $\Xi_c$ is chosen as product of one-dimensional interpolation nodes (such as Chebyshev nodes \cite{FONG20098712}). Hence, for any source cell $s$, the IBFMM M2M operator between $s$ and $s'$, with $s'$ a daughter of $s$, is defined as
\begin{equation}
\label{eq:M2M}
    V[s,s'] := \begin{bmatrix}S_0[s](\bo{y}^{s'}_0)& \hdots & S_0[s](\bo{y}^{s'}_{P-1})\\
    \vdots&\ddots&\vdots\\
    S_{P-1}[s](\bo{y}^{s'}_0)& \hdots & S_{P-1}[s](\bo{y}^{s'}_{P-1})\end{bmatrix}.
\end{equation}

Using Eq. \eqref{interpolationformula}, we can identify that the IBFMM L2P operator on leaf cell $c$ is the adjoint of the IBFMM P2M operator on this $c$. Hence, for any non-leaf cell $c$, the IBFMM L2L operator between $c$ and a daughter $c'$ of $c$ is the adjoint of the IBFMM M2M operator between $c$ and $c'$, that is:
\begin{equation}
\label{eq:L2P_L2L}
    U^*[c] = V[c]^*,\hspace{0.3cm}U^*[c',c] = V[c,c']^*.
\end{equation}

\begin{figure}
    \centering
    \includegraphics[width=\linewidth]{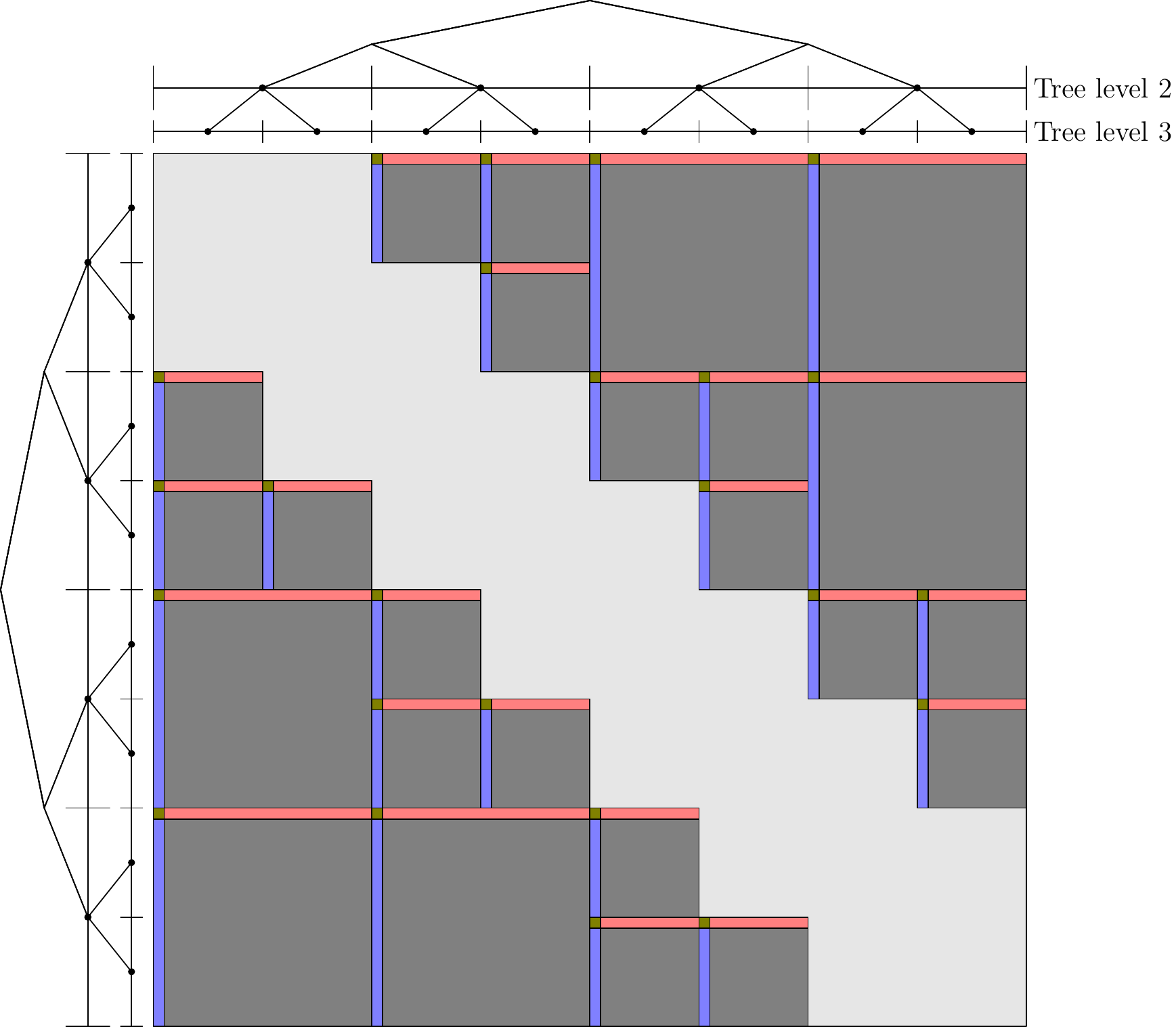}
    \caption{Schematic representation of low-rank approximations obtained through interpolation over 1D quasi-uniform particle distributions. Perfect trees (binary trees here) are represented above (sources) and on the left (targets) of the matrix representation of the far field compression. There are admissible cell pairs $(t,s)$ only at levels 2 (big dark grey blocks) and 3 (small dark grey blacks). In each of the corresponding blocks, we drew its low-rank approxiation using red color for source polynomials (i.e. $V[s]$ matrix), blue color for target polynomials (i.e. $U^*[t]$ matrix and purple for small square matrices $C[t,s]$). Non-admissible cell pairs are represented using light grey.\label{fig:ibfmmin1d}}
\end{figure}

As for the classical FMM algorithm, the IBFMM reaches the linear complexity with respect to the number of particles in the system. In the remainder of this article, we consider that $K$ is translationally and rotationally invariant, which is always verified in our applications. A schematic view of IBFMM is provided in Fig. \ref{fig:ibfmmin1d}.

\section{ANKH-FMM}
\label{s::ankhfmm}
Numerical issues arise when dealing with Eq. \ref{eq_qmutheta_coulomb} since the series only conditionally converges, possibly preventing practical convergence. The Ewald summation technique \cite{pmeref} transforms Eq. \eqref{eq_qmutheta_coulomb} into $\cc{E} = \cc{E}_{real} + \cc{E}_{rec} + \cc{E}_{self}$
where
\begin{equation}
\label{eq:ewald}
\begin{aligned}
    \cc{E}_{real} :=& \sum_{\bo{t}}\sum_{\bo{x}}\sum_{\bo{y}}\kk{D}_\bo{x}\kk{D}_\bo{y}\left(\frac{erfc\left(\xi|\bo{x}-\bo{y}+\bo{t}|\right)}{|\bo{x}-\bo{y}+\bo{t}|}\right),\\
    \cc{E}_{rec} :=& \frac{1}{2\pi |B|(2r_B)^2} \sum_{\bo{m}\neq \bo{0}}\frac{e^{-2\left(\frac{\pi r_B}{\xi}\right)^2\bo{m}\cdot \bo{m}}}{\bo{m}\cdot \bo{m}}S(\bo{m})S(-\bo{m}),\\
S(\bo{m}) :=& \sum_{\bo{x}}e^{i\frac{\pi}{r_B}\bo{m}\cdot \bo{x}}\Big(q_\bo{x} + 2i\pi \hspace{0.07cm}r_B\hspace{0.05cm}\mu_\bo{x}\cdot \bo{m}\\& - \big(2\pi \hspace{0.07cm}r_B)\big)^2\hspace{0.05cm}\Theta_{\bo{x}}:\left(\bo{m}\bo{m}^T\right)\Big),\\
    \cc{E}_{self} :=& -\frac{\xi}{\sqrt{\pi}}\sum_{\bo{x}}\Big(q_\bo{x}^2 + \frac{2\xi^2}{3}\mu_\bo{x}\cdot \mu_\bo{x} + \frac{8\xi^4}{5}\Theta_\bo{x}:\Theta_\bo{x}\Big),
\end{aligned}
\end{equation}
in which all the series are absolutely convergent, $\xi\in \bb{R}^{+*}$ and $erfc(x) = 1-\frac{2}{\pi}\int_{0}^{x}e^{-t^2}dt$ denotes the complementary error function. The main point is that $erfc$ quickly decays for sufficiently large $\xi$, implying that only a small amount of images (i.e. of different $\bo{t}$'s) in $\cc{E}_{real}$ have to be considered to reach a given accuracy. Although, applying a cutoff with radius $\delta$ on the sum $\sum_\bo{x}\sum_\bo{y}$ transforms it into $\sum_\bo{x}\sum_{\bo{y}\in \cc{B}(\bo{x},\delta)}$ accelerating its evaluation, using the notation $\cc{B}(\bo{x},\delta) := \{\bo{z}\in \bb{R}^3\hspace{0.02cm}|\hspace{0.02cm}|\bo{z}-\bo{x}|\leq \delta\}$. $\cc{E}_{rec}$ also is a quickly converging sum and only a limited number of $\bo{m}$'s have to be considered to reach a given accuracy, provided that $\xi$ is sufficiently large. Usually, $\xi$ is chosen to minimize the evaluation cost according to the possible cutoffs, balancing the computation costs of $\cc{E}_{real}$ and $\cc{E}_{rec}$. The Particle Mesh Ewald \cite{pmeref} (PME) algorithm uses Fast Fourier Transforms (FFTs) in the computation of $S(\bo{m})$'s, further reducing the complexity.

\subsection{A new alternative to PME}
\label{subsect:newalternative}
On one side, the parallelization of PME algorithm in distributed memory can be a practical bottleneck due to the scalability of the FFT \cite{10.1007/978-3-030-86359-3_21}. On the other side, there exists fast scalable methods in distributed memory for the evaluation of $N$-body problems involving asymptotically smooth kernels, such as (IB)FMM (see Sect. \ref{ss::ibfmmkis}). We thus propose a new alternative to PME, based on a very simple idea. Indeed, the kernel
\begin{equation}
\label{eq:Hdef}
    H(\bo{x},\bo{y}) := \frac{erfc\left(\xi|\bo{x}-\bo{y}|\right)}{|\bo{x}-\bo{y}|},
\end{equation}
is asymptotically smooth \cite{badreddine:hal-03826012}, as well as all its derivatives. For sufficiently large $p$, thanks to the absolute convergence of the terms in Eq. \eqref{eq:ewald}, we have
\begin{equation}
\label{eq_real_nbody}
\begin{aligned}
    \cc{E}_{real} &\approx \sum_{\bo{t}\in 2r_B Z_p}\sum_{\bo{x}}\sum_{\bo{y}}\kk{D}_\bo{x}\kk{D}_\bo{y}\left(\frac{erfc\left(\xi|\bo{x}-\bo{y}+\bo{t}|\right)}{|\bo{x}-\bo{y}+\bo{t}|}\right)\\
    &= \sum_{\bo{t}\in 2r_B Z_p}\sum_{\bo{x}}\sum_{\bo{y}}\kk{D}_\bo{x}\kk{D}_\bo{y}H(\bo{x},\bo{y}+\bo{t}) =: \cc{\tilde{E}}_{real},
\end{aligned}
\end{equation}
with $Z_p := \big\{(z_0,z_1,z_2)\in \bb{Z}^3\hspace{0.1cm}|\hspace{0.1cm}|z_k|\leq p,\hspace{0.1cm}k\in [\![0,2]\!]\big\}$ and $p$ to be fixed by user.

Our approach is thus very simple: we consider small $\xi$ (typically $\xi = 0.01$), so that $\cc{E}_{rec}\approx 0$ and can be neglected for our targeted accuracies. In this case, the counterpart relies in the computation of $\cc{\tilde{E}}_{real}$ that requires a relatively large amount of periodic images (i.e. a large $p$) to reach the numerical convergence. Since $\cc{E}_{self}$ involves $\cc{O}(N)$ flops to be computed (and is embarrassingly parallel), only this $\cc{\tilde{E}}_{real}$ has a corresponding intensive computational cost. Hopefully, Eq. \eqref{eq_real_nbody} allows the computation of $\tilde{\cc{E}}_{real}$ to be handled through fast methods for $N$-body problems with kernel $\kk{D}_\bo{x}\kk{D}_\bo{y}H(\bo{x},\bo{y})$.

\subsection{Reformulation}
\label{ss:reformulation}
In this section, we consider that $p=0$. Extension to other $p$'s will be done in Sect. \ref{ss:fmmperiodicity}. Let $\cc{\tilde{E}}_{real}(t,s):=\sum_{\bo{x}\in t}\sum_{\bo{y}\in s}\kk{D}_\bo{x}\kk{D}_\bo{y}H(\bo{x},\bo{y})$, 
where $t,s\in \cc{L}(\scr{T})$. We can now divide $\cc{\tilde{E}}_{real}$ into two parts:
\begin{equation*}
    \underbrace{\sum_{\substack{t,s\in \cc{L}(\scr{T})\\ dist(t,s)= 0}}\cc{\tilde{E}}_{real}(t,s)}_{=: \cc{N}_{real}} + \underbrace{\sum_{\substack{t,s\in \cc{L}(\scr{T})\\ dist(t,s)\neq 0}}\cc{\tilde{E}}_{real}(t,s)}_{=: \cc{F}_{real}}
\end{equation*}
where $dist(t,s)$ denotes the distance between $t$ and $s$. Hence, by increasing the octree depth, one decreases the number of interactions computed inside $\cc{N}_{real}$ (\textit{near field}) while increasing the number of interactions in $\cc{F}_{real}$ (\textit{far field}).

Each $(t,s)\in \cc{L}(\scr{T})$ such that $dist(t,s)\neq 0$ actually is a well-separated pair of cells. Hence, by means of multivariate polynomial interpolation as in Eq. \eqref{interpolationformula} on $t$ and $s$, one may accurately approximate $\cc{\tilde{E}}_{real}(t,s)$ as: 
\begin{equation}
\label{eq:M_t}
    \begin{aligned}
         & \sum_{\bo{x}\in t}\sum_{\bo{y}\in s}\kk{D}_\bo{x}\kk{D}_\bo{y}\left(\sum_{k}S_k[t](\bo{x})\sum_{l}H(\bo{x}_k^t,\bo{y}_l^s)S_l[s](\bo{y})\right)\\
        =& \sum_{k}\underbrace{\left(\sum_{\bo{x}\in t}\kk{D}_\bo{x}S_k[t](\bo{x})\right)}_{=:\cc{Q}_t(\bo{x}_k^t)}\sum_{l}H(\bo{x}_k^t,\bo{y}_l^s)\underbrace{\left(\sum_{\bo{y}\in s}\kk{D}_\bo{y}S_l[s](\bo{y})\right)}_{=:\cc{Q}_s(\bo{y}_l^s)}\\
        =&\sum_{k}\cc{Q}_t(\bo{x}_k^t)\sum_{l}H(\bo{x}_k^t,\bo{y}_l^s)\cc{Q}_s(\bo{y}_l^s).
    \end{aligned}
\end{equation}
These $\cc{Q}_c(\bo{u})$, with $\bo{u}$ any interpolation node in $c$, can be interpreted as modified charges associated to this interpolation node, which can itself be seen as a particle. Hence, one can perform an IBFMM in which the multipole/local expansions are computed using the modified charges while the near field (i.e. $\cc{N}_{real}$) is evaluated using the explicit formula for the complete interaction between particles with respect to the kernel $H$ (i.e. involving charges, dipoles and quadrupoles here). In other words, this method (that we named ANKH) consists in the following steps:
\begin{enumerate}
    \item Compute the tree decomposition $\scr{T}$ of the simulation box $B$,
    \item For each leaf $c\in \cc{L}(\scr{T})$, compute the modified charges in $c$ (using the same interpolation grid up to a translation to the center of $c$),
    \item Compute $\cc{N}_{real}$ using direct computation,
    \item Compute $\cc{F}_{real}$ solving the $N$-body problem defined on interpolation nodes and modified charges,
    \item Compute $\cc{E}_{self}$ using Eq. \eqref{eq:ewald},
    \item Return $\cc{N}_{real} + \cc{F}_{real} + \cc{E}_{self}$.
\end{enumerate}
The $N$-body problem mentioned in step 4 can be written this way:
\begin{equation}
\label{eq:reformulationnnbody}
    \sum_{\substack{t\in \cc{L}(\scr{T}),\\\bo{u}\text{ interpolation}\\\text{node on }t}}\cc{Q}_t(\bo{u})\underbrace{\left(\sum_{\substack{s\in \cc{L}(\scr{T}),\\\bo{v}\text{ interpolation}\\\text{node on }s}}H_{t,s}(\bo{u},\bo{v})\cc{Q}_s(\bo{v})\right)}_{N\text{-body problem as in Eq. \eqref{nbodypb}}}
\end{equation}
where $H_{t,s} = \begin{cases}H\text{ if }dist(t,s)> 0,\\0\text{ otherwise}\end{cases}$. According to Eq. \eqref{eq:reformulationnnbody}, the quantity into parenthesise can be efficiently computed in linear time by means of IBFMM applied on kernel $H$ and in which we simply avoid the P2P operators. We name this approach ANKH-FMM.

In the case of quasi-uniform particle distribution (actually the kind of particle distributions we are targeting), a perfect octree with depth $~log(N)$ generates leaves with $\cc{O}(1)$ particles. Hence, in this situation, step 3 of the ANKH method costs $\cc{O}(N)$ flops. For such perfect trees, the tree construction (step 1) also is $\cc{O}(N)$ and computation of $\cc{E}_{self}$ is, of course, also linear in complexity (step 5). Since it is obvious that step 6 is $\cc{O}(1)$, the overall complexity is $\cc{O}(N)$.

Two lasts questions remain: what about the periodicity and the computation of modified charges? The first of these questions is the purpose of Sect. \ref{ss:fmmperiodicity} and the second one is the purpose of Sect. \ref{ss:numericaldiff}.

\subsection{Efficient handling of periodicity}
\label{ss:fmmperiodicity}
Various formulations of the FMM provide efficient ways of dealing with periodicity \cite{Schmidt1991ImplementingTF,challacombe,407723}. We consider a slightly different approach in our method, well-suited for IBFMM. One idea is to perform a classical IBFMM on the computational box as well as on the direct $3^3-1$ adjacent images of it. Then, all other distant images are well-separated from the computational box and their interaction can be approximated using only the multipole expansion $\cc{M}_B$ of this box. More precisely, denoting $C[B,I]$ the M2L matrix between the computational box $B$ and its image $I$, one wants to fastly compute:
\begin{equation}
\label{eq:sumimages}
    \sum_{\substack{I\in Z_p\\p\geq |I|>1}}C[B,I]\cc{M}_B = \left(\sum_{\substack{I\in Z_p\\p\geq |I|>1}}C[B,I]\right)\cc{M}_B.
\end{equation}
The remaining task is thus to sum the M2L matrices between images of B and B itself efficiently. This can be done by means of interpolation over equispaced grid \cite{blanchard:tel-01534930,chollet:hal-03563005}.

Indeed, such tool allow to express $C[B,I]$ as a product $\chi^*\bb{F}^*D[B,I]\bb{F}\chi$, where $\chi$ is a zero-padding matrix that increases the size of the inner matrix by slightly less than $8$, $\bb{F}$ denotes the discrete Fourier matrix (that can be applied through FFT) and $D[B,I]$ is a diagonal matrix that can be computed in linearithmic time with respect to the size of $C[B,I]$.

More precisely, let $\Xi_{B}:=\{\Xi_{B}(i,j,k):=(-r_B+2r_Bi/(L-1),-r_B+2r_Bj/(L-1),-r_B+2r_Bk/(L-1))\hspace{0.01cm}|\hspace{0.01cm}i,j,k\in [\![0,L-1]\!]\}$ be a equispaced grid over $B$ with $L$ nodes in each direction (i.e. with $L^3$ nodes in total) and let $\Xi_I$ be its translation from the origin to the center of the image $I$ of $B$. The matrix $C[B,I]$, in this case, has a block-Toeplitz structure \cite{blanchard:tel-01534930,chollet:hal-03563005}. Any such block-Toeplitz matrix can be embedded into a circulant one of size $(2L-1)^3$ whose first row is composed of the $(2L-1)^3$ different entries of $C[B,I]$. $\chi\in\{0,1\}^{(2L-1)^3\times L^3}$ and $D[B,I]\in \bb{C}^{(2L-1)^3\times (2L-1)^3}$ are then defined as follows:
\begin{equation}
\label{eq:circdef}
\begin{aligned}
    &\chi_{q,p} := \begin{cases}1\text{ if }\substack{i_0=i_1,j_0=j_1;\\k_0=k_1;i_1,j_1,k_1<L}\\0\text{ otherwise}\end{cases},\\
    &D[B,I] := diag(\bb{F}R[B,I]),
\end{aligned}
\end{equation}
where $R[B,I]$ is the first row of the circulant matrix formed by embedding of the Toeplitz matrix $C[B,I]$, i.e.
\begin{equation}
\label{eq:defRz}
\begin{aligned}
    &R[B,I]_{i(2L-1)^2+j(2L-1)+k}:=H\left(\bo{z}^{(I)}_{i,j,k},\bo{0}\right),\\
    &\bo{z}^{(I)}_{ijk}:=\Xi_{B}(\tilde{i},\tilde{j},\tilde{k}) - \Xi_I(\hat{i},\hat{j},\hat{k}),\\
    &\tilde{i} := \begin{cases}0\text{ if }i<L\\2L-1-i\text{ otherwise}\end{cases}\\
    &\hat{i} := \begin{cases}i\text{ if }i<L\\0\text{ otherwise}\end{cases}
\end{aligned}
\end{equation}
Since $\bb{F}$ is practically performed through FFT, computation of $D[B,I]$ results in $\cc{O}(L^3log(L))$ flops.

We can thus reformulate Eq. \eqref{eq:sumimages} as
\begin{equation}
\label{eq:nsumimages}
    \chi^*\bb{F}^*\underbrace{\left(\sum_{\substack{I\in Z_p\\p\geq |I|>1}}D[B,I]\right)}_{=:D_B}\bb{F}\chi\cc{M}_B.
\end{equation}
As a sum of diagonal matrix, $D_B$ clearly is diagonal. We can switch the sum in Eq. \eqref{eq:nsumimages} to the vector $R[B,I]$ thanks to Eq. \eqref{eq:circdef}, resulting into
\begin{equation}
    D_B=diag\left(\bb{F}\sum_{\substack{I\in Z_p\\p\geq |I|>1}}R[B,I]\right).
\end{equation}

\begin{figure}
    \centering
    \includegraphics[width=\linewidth]{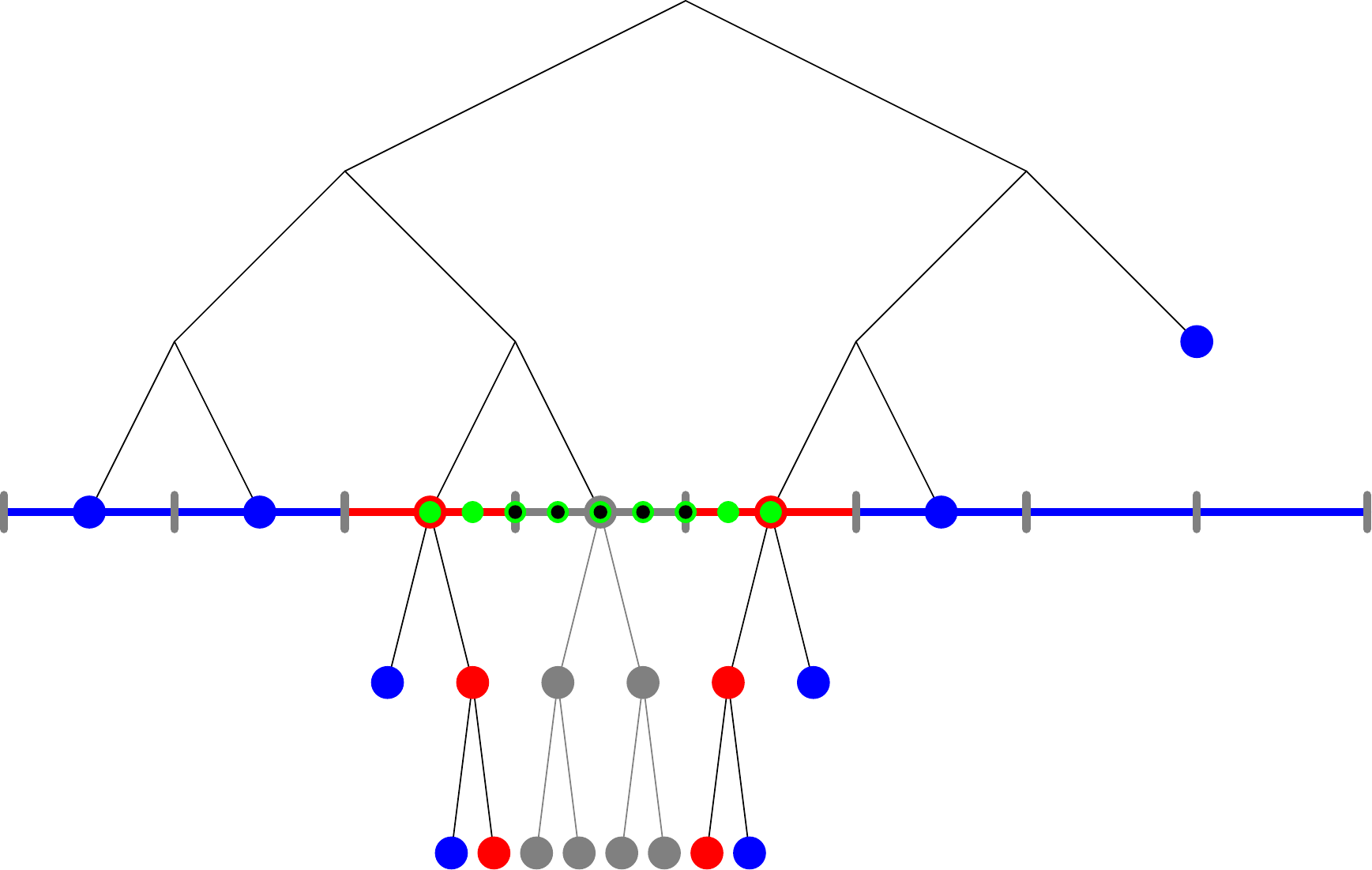}
    \caption{1D schematic representation of the handling of far periodic images through hierarchical methods. Octree $\scr{T}$ in the simulation box $B$ is represented in grey. Its two direct (non well-separated) neighbor images (in red) directly interact with $\scr{T}$. Red and blue nodes in these trees respectively correspond to non well-separated and well-separated nodes of the images regarding nodes of $\scr{T}$ at the same tree level. On the root of $\scr{T}$ (i.e. on $B$), interpolation nodes over a equispaced grid of size 5 are represented (in black). $P_{target}$ is drawn in green. Blue nodes at the same tree level than $B$ are elements of $P_{source}$. Upper nodes can be interpreted as aggregation of some of them during a hierarchical procedure to solve the $N$-body problem of Eq. \eqref{eq::nbodyperioima}.\label{fig:periodic_tree}}
\end{figure}

Since the image $I$ of $B$ is a translation of $B$ by $\bo{t}=2r_BI$, we have that $\Xi_I(\hat{i},\hat{j},\hat{k}) = \Xi_B(\hat{i},\hat{j},\hat{k})+\bo{t}$, so that, thanks to Eq. \eqref{eq:defRz} and the translational invariance of $H$,
\begin{equation}
\label{eq:prenbodyforimages}
\begin{aligned}
    R[B,I]_{i(2L-1)^2+j(2L-1)+k}&=H\left(\bo{z}^{(\bo{0})}_{i,j,k}+\bo{t},\bo{0}\right)\\
    &=H\left(\bo{z}^{(\bo{0})}_{i,j,k},-\bo{t}\right).
\end{aligned}
\end{equation}
Notice that $-I$ also lies into the set of images of $B$ so that one can eliminate the minus sign in Eq. \eqref{eq:prenbodyforimages}. Also, $\{\bo{z}^{(\bo{0})}_{i,j,k}\hspace{0.01cm}|\hspace{0.01cm}i,j,k\in [\![0,(2L-1)]\!]\}$ can be interpreted as a particle distribution over $B$ and its direct neighbors. Moreover, applying the sum over $I$ of Eq. \eqref{eq:nsumimages} to Eq. \eqref{eq:prenbodyforimages} gives
\begin{equation}
\label{eq::nbodyperioima}
    \sum_{\substack{I\in Z_p\\p\geq |I|>1}}H\left(\bo{z}^{(\bo{0})}_{i,j,k},\bo{t}\right)
\end{equation}
which can be interpreted as a $N$-body problem with all charges set to $1$ over two particle distributions: the set $P_{target}$ of $(2L-1)^3$ different $\bo{z}^{(\bo{0})}_{ijk}$'s and the set $P_{source}$ of $\bo{t}$'s. For large $p$, using a hierarchical method for such problem, one may compute this quantities very efficiently since the first of these particle distribution is already well-separated from each $\bo{t}$. In addition, both $H$, $P_{target}$ and $P_{source}$ are invariant under the action of the rotation group that preserves the cube. This can be exploited to speed-up the computation \cite{chollet:hal-03563005}. See Fig. \ref{fig:periodic_tree} for a schematic representation.

A critical point here (that actually motivates this entire subsection) is that $D_B$ does not depend on the particle distribution nor its modified charges, but only relies on the size of the simulation box $B$. Hence, the computation of $D_B$ can be done only once inside a precomputation step, not at each IBFMM call. This has a particular interest since in daylife molecular dynamics applications, such IBFMM calls on the same $B$ would be numerous, but needing only the single precomputation of $D_B$ with our method.

In terms of complexity, at runtime, $D[B,I]$ is already computed, so the cost of handling the periodicity this way is the one of two FFTs on a small grid of size $(2L-1)^3$. Since $L$ is a (small) constant, these FFTs cost $\cc{O}(L^3log(L)) = \cc{O}(1)$ flops.

\subsection{Numerical differentiation}
\label{ss:numericaldiff}
The idea of switching derivatives to the interpolation polynomials has already found applications in the literature \cite{messnerthesis}. In our case, things are slightly more difficult because we consider differential operators instead of simple derivatives. To do so, we consider a Lagrange interpolation over product of 1D Chebyshev grids in each leaf cell as well as all along the tree structure with a fixed order $L$ in each direction. Notice that this interpolation is different from the one of Sect. \ref{ss:fmmperiodicity}, that performs over a equispaced grid and only for the root. However, this does not impact the method at all since the root is only involved in interactions between non-adjacent images of $B$ (see Fig. \ref{fig:periodic_tree}) and because the definition of the M2M/L2L operators is valid no matter the interpolation grid type is.

The choice of multivariate Lagrange interpolation over products of 1D Chebyshev nodes is motivated by the numerical instability that arises on large interpolation order and at fixed arithmetic precision when dealing with equispaced grids. Moreover, this choice of interpolation is known to be particularly good in the context of multilevel summation schemes.

As presented in Sect. \ref{ss:reformulation}, in the ANKH approach, the differential operators apply on interpolation polynomials directly. Any Lagrange polynomial over a product grid is the product of the Lagrange polynomials in each one-dimensional grids. Assuming for the sake of clarity (One only has to introduce a scaling function in order to recover the general case \cite{messnerthesis,blanchard:tel-01534930}) that the one-dimensional grids are defined on $[-1,1]$, we thus have the following decomposition of such polynomial $S$:
\begin{equation}
\label{eq:multivariatepolynomials}
    S_{iL^2+jL+k}(\bo{x}) = S_i(x_0)S_j(x_1)S_k(x_1)
\end{equation}
where $S_{iL^2+jL+k}$ can be either $S_k[t]$ or $S_l[s]$ following the notations of Eq. \eqref{interpolationformula}, $S_i$'s are one-dimensional Lagrange polynomials over Chebyshev (i.e. Lagrange-Chebyshev polynomials) nodes and $\bo{x}=(x_0,x_1,x_2)\in \bb{R}^3$. Explicit forms \cite{FONG20098712} are known for these polynomials $S_i$:
\begin{equation}
\label{eq:lagrangechebinterp}
    S_k(x) = \frac{1}{L} + \frac{2}{L}\sum_{m=1}^{L-1}T_m(r_k)T_m(x)
\end{equation}
where $T_m$ is the $m^{th}$ Chebyshev polynomial and $r_k := cos\left(\frac{2k-1}{2L}\right)$ is the $k^{th}$ Chebyshev node. $T_m(x)$ can be computed using the recurrence relation:
\begin{equation}
\label{eq:recformulacheby}
\begin{aligned}
    &T_{m+2}(x) = 2xT_{m+1}(x) - T_{m}(x), \\&T_0(x)=1, T_1(x)=x.
\end{aligned}
\end{equation}

Problem is that explicit derivation of Eq. \eqref{eq:lagrangechebinterp} leads to 
\begin{equation}
        S_k'(x) \frac{2}{L\sqrt{1-x^2}}\sum_{m=1}^{L-1}mT_m(x_k)sin(m\hspace{0.1cm}acos(x))
\end{equation}
which is numerically unstable near the bounds of $[-1,1]$. This can here be solved by using Chebyshev polynomials of the second kind in the derivation, but the same kind of problem appears when trying to derive $S_k$ once again (which is needed to take into account quadrupoles). We thus opted for another way of deriving these polynomials: there exists \cite{Amiraslani2019DifferentiationMF} a upper triangular matrix $\cc{D}$ such that:
\begin{equation}
\label{eq:derivemanytimes}
    \bo{T}'(x) = \cc{D}\bo{T}(x), \hspace{0.2cm}\bo{T}(x):=\begin{bmatrix}T_0(x)\\\vdots\\T_{L-1}(x)\end{bmatrix}.
\end{equation}
This can trivially be generalized to any derivation order, giving:
\begin{equation}
\label{eq:generalderivationcheby}
    \bo{T}^{(n)}(x) = \cc{D}^n\bo{T}(x).
\end{equation}
This last form does not suffer from numerical instability near the bounds of $[-1,1]$. To exploit it in ANKH, we start by reformulating the set of Lagrange-Chebyshev polynomials of Eq. \eqref{eq:lagrangechebinterp} in matrix form
\begin{equation}
\label{eq:lagrangechebymatrix}
\begin{aligned}
    &\bo{S}(x) = \bb{H}\bo{T}(x),\\&\bb{H} := \frac{1}{L}\begin{bmatrix}
      1 & 2T_1(r_0) & \hdots & 2T_{L-1}(r_0)\\
      \vdots & \vdots & \ddots & \vdots\\
      1 & 2T_1(r_{L-1}) & \hdots & 2T_{L-1}(r_{L-1})\\
    \end{bmatrix},
\end{aligned}
\end{equation}
where $\bb{H}$ does not depend on $x$ and be be precomputed only once. We can now combine Eq. \eqref{eq:generalderivationcheby} and Eq. \eqref{eq:lagrangechebymatrix} in order to obtain our general formula:
\begin{equation}
    \bo{S}^{(n)}(x) = \bb{H}\cc{D}^n\bo{T}(x).
\end{equation}
To extend this technique to the three-dimensional case, a direct application of Eq. \eqref{eq:multivariatepolynomials} provides expressions of the multivariate polynomials.

We summarize using $C$++-like pseudocode in Alg. \ref{alg:compute_der_p_S} the steps needed to compute the modified charges for multipoles up to quadrupoles (trivial generalizations can be extrapolated from this pseudo-code) using numerical differentiation through derivation of Chebyshev polynomials. Notice that the different $\bo{T}$'s computed in Alg. \ref{alg:compute_der_p_S} could be stacked in practice inside a matrix with $9$ columns before applying the product by $\bb{H}$, which allows to benefit from matrix-matrix products instead of matrix-vector ones (i.e. BLAS3 instead of BLAS2), that should better perform. 

\begin{algorithm}[htbp]
\caption{Compute $\cc{Q}_s$ in leaf cell $s$}
\label{alg:compute_der_p_S}
\begin{algorithmic}
\State // Initialize modified charges with zero
\State $\cc{Q}_s = 0;$
\State // Loop over particles in $s$ (order does not matter)
\For{$\bo{y}\in Y_{|s}$}
\State // Notation: $\bo{y} = (y_0,y_1,y_2)$
\State // Allocate arrays storing Chebyshev polynomials
\State $S_0^{(0)},S_0^{(1)},S_0^{(2)} = \text{new }\bb{R}[L]$;
\State $S_1^{(0)},S_1^{(1)},S_1^{(2)} = \text{new }\bb{R}[L]$;
\State $S_2^{(0)},S_2^{(1)},S_2^{(2)} = \text{new }\bb{R}[L]$;
\State $\bo{T} \hspace{1.67cm}=  \text{new }\bb{R}[L];$
\State // Decompose the computation per direction
\For{$k\in[\![0,2]\!]$}\Comment{$\cc{O}(L^2)$ flops}
    \State // Initialize the first two polynomials, see Eq. \eqref{eq:recformulacheby}
    \State $\bo{T}[0] = 1$;
    \State $\bo{T}[1] = y_k$;
    \State // Compute recurrence following Eq. \eqref{eq:recformulacheby}
    \For{int $i=2$; $i<L$; $i++$} \Comment{$\cc{O}(L)$ flops}
        \State $T[i] = 2y_kT_k[i-1]-T_k[i-2]$;\Comment{$\cc{O}(1)$ flops}
    \EndFor
    \State // Apply differentiation matrix following Eq. \eqref{eq:derivemanytimes}
    \For{int $j=0$; $j<3$; $j++$} \Comment{$\cc{O}(L^2)$ flops}
        \State $\bo{T} \hspace{0.37cm}= \cc{D}\cdot \bo{T}$; \Comment{$\cc{O}(L^2)$}
        \State $S_k^{(p)} = \bb{H}\cdot \bo{T}$; \Comment{$\cc{O}(L^2)$ flops}
    \EndFor
\EndFor
\State // Add all possible combinations scaled with charge...
\State $\cc{Q}_s = \cc{Q}_s + S^{(0)}_0S^{(0)}_1S^{(0)}_2q_\bo{y};$
\State // ... dipoles...
\State $\cc{Q}_s = \cc{Q}_s + S^{(1)}_0S^{(0)}_1S^{(0)}_2\left(\mu_\bo{y}\right)_0;$
\State $\cc{Q}_s = \cc{Q}_s + S^{(0)}_0S^{(1)}_1S^{(0)}_2\left(\mu_\bo{y}\right)_1;$
\State $\cc{Q}_s = \cc{Q}_s + S^{(0)}_0S^{(0)}_1S^{(1)}_2\left(\mu_\bo{y}\right)_2;$
\State // ... and quadrupoles, exploiting symmetry
\State $\cc{Q}_s = \cc{Q}_s + S^{(2)}_0S^{(0)}_1S^{(0)}_2\left(\Theta_\bo{y}\right)_{0,0};$
\State $\cc{Q}_s = \cc{Q}_s + S^{(1)}_0S^{(1)}_1S^{(0)}_2\left(\Theta_\bo{y}\right)_{0,1}\times 2;$
\State $\cc{Q}_s = \cc{Q}_s + S^{(1)}_0S^{(0)}_1S^{(1)}_2\left(\Theta_\bo{y}\right)_{0,2}\times 2;$
\State $\cc{Q}_s = \cc{Q}_s + S^{(0)}_0S^{(2)}_1S^{(0)}_2\left(\Theta_\bo{y}\right)_{1,1};$
\State $\cc{Q}_s = \cc{Q}_s + S^{(0)}_0S^{(1)}_1S^{(1)}_2\left(\Theta_\bo{y}\right)_{1,2}\times 2;$
\State $\cc{Q}_s = \cc{Q}_s + S^{(0)}_0S^{(0)}_1S^{(2)}_2\left(\Theta_\bo{y}\right)_{2,2};$

\EndFor
\end{algorithmic}
\end{algorithm}

Regarding the total complexity, the computation of $Q$ inside a leaf cell $s$ has a cost of $\cc{O}\left(L^3N_s\right)$, where $N_s$ denotes the number of particles in $s$. Since $L$ is a constant, one may drop it from the big $\cc{O}$ and the total complexity becomes:
\begin{equation}
    \sum_{s\in \cc{L}(\scr{T})}\cc{O}\left(N_s\right) = \cc{O}\left(\sum_{s\in \cc{L}(\scr{T})}N_s\right) = \cc{O}(N).
\end{equation}

\subsection{Overall complexity and advantages}
\label{ss:complexity_ankh_fmm}
At the end of the day, ANKH-FMM defines a linear method since the overall complexity, thanks to all previous sections, can be counted this way:
\begin{equation}
    \underbrace{\cc{O}(N)}_{\text{step }1}+\underbrace{\cc{O}(N)}_{\text{step }2}+\underbrace{\cc{O}(N)}_{\text{step }3}+\underbrace{\cc{O}(N)}_{\text{step }4}+\underbrace{\cc{O}(N)}_{\text{step }5}+\underbrace{\cc{O}(1)}_{\text{step }6}
\end{equation}
where the steps are those presented in Sect. \ref{ss:reformulation}. This has to be compared with the linearithmic complexity of PME. Hence, our method is asymptotically less costly than PME. However, the FMM (so as IBFMM) is known to suffer from an important prefactor, while the FFT has a (very) small one and may almost reach CPU peak performance. This means that the theoretical complexities are not sufficient to compare the two approaches and numerical tests (that will be provided in the following) have to verify the efficiency of ANKH.

However, as we mentioned in Sect. \ref{subsect:newalternative}, there is another strong advantage of using ANKH-FMM: its distributed parallel potential. Indeed, computation of modified charges, as local operations, are easy to parallelize in distributed memory context. Moreover, IBFMM has shown impressive parallel performance in this context and general parallel strategies aiming at reaching exascale are already developed for FMMs. Parallel scaling of FMMs can be put in comparison with ones of FFTs on which PME relies. In addition, efficient implementation of IBFMM on GPU \cite{gpudarve} were proposed in the literature for perfect trees (i.e. exactly our application case). Hence, the ANKH strategy appears as a serious theoretical candidate for exascale hybrid computations.

We may already mention that GPU implementations of IBFMM are technical. This is why we wanted to devise a simple and portable scheme to run on GPU while keeping the parallel structure of the FMM. This is the purpose of Sect. \ref{s:ankhfft}.

\subsection{Mutual interactions}
\label{ss:mutual}
In the particular context of energy computation, the global algorithmic discussed in Sect. \ref{ss:complexity_ankh_fmm} still is suboptimal. Indeed, coming back to Eq. \eqref{eq:fmmchain}, we locally want to compute $q_{|t'}^T\phi_{t',s'}$ instead of $\phi_{t'}$ only (that would rather be suited for forces computations). Hence, we would have
\begin{equation}
\begin{aligned}
    q_{|t'}^T\phi_{t',s'} &\approx \left(q_{|t'}^T U^*[t']\right)U^*[t',t]C[t,s]V[s,s']\underbrace{\left(V[s']q_{|s'}\right)}_{=\cc{M}_{s'}}\\
    &=\left(q_{|t'}^T V[t']^*\right)V[t,t']^*C[t,s]V[s,s']\cc{M}_{s'}\\
    &= \left(V[t']q_{|t'}\right)^TV[t,t']^*C[t,s]V[s,s']\cc{M}_s\\
    &= \cc{M}_{t'}^TV[t',t]^*C[t,s]V[s,s']\cc{M}_{s'}\\
    &= \underbrace{\left(\cc{M}_{t'}^TV[t,t']^T\right)}_{=\left(V[t,t']\cc{M}_t\right)^T = \cc{M}_t^T}C[t,s]\underbrace{\left(V[s,s']\cc{M}_{s'}\right)}_{=\cc{M}_s}
\end{aligned}
\end{equation}
since $U^*[t']=V[t']^*=V[t']^T$ because IBFMM L2P and P2M operators are adjoint matrices with real entries (so as IBFMM L2L and M2M operators). This provides a new formula for $\cc{F}_{real}$:
\begin{equation}
\label{eq:refomulatefreal}
\begin{aligned}
    \cc{F}_{real} &\approx \sum_{\text{target cell }t}\cc{M}_t^T\sum_{s\in \Lambda(s)}C[t,s]\cc{M}_s\\
    &=\sum_{\text{target cell }t}\sum_{s\in \Lambda(s)}\left(\cc{M}_t^TC[t,s]\cc{M}_s\right).
\end{aligned}
\end{equation}
This does not reduce the theoretical complexity but strongly simplifies the method. In addition, many optimizations may be derived from this formula.

First, when performing M2L, we can directly compute the quantity into parenthesis in Eq. \eqref{eq:refomulatefreal}, i.e. adding the left product by $\cc{M}_t^T$. Second, since if two cells $t,s\subseteq B$ are such that $s\in \Lambda(t)$, then $t\in \Lambda(s)$, meaning that the two interactions will be computed. This can be simplified because $C[t,s] = C[s,t]^T$:
\begin{equation}
\label{eq;mutualM2L}
\begin{aligned}
    &\cc{M}_t^TC[t,s]\cc{M}_s + \cc{M}_s^TC[s,t]\cc{M}_t\\
    =&\cc{M}_t^TC[s,t]^T\cc{M}_s + \cc{M}_s^TC[s,t]\cc{M}_t\\
    =&\left(\cc{M}_s^TC[s,t]\cc{M}_t\right)^T + \cc{M}_s^TC[s,t]\cc{M}_t\\
    =& \hspace{0.1cm}2\cc{M}_s^TC[s,t]\cc{M}_t.
\end{aligned}
\end{equation}
This implies that only half of the M2L interactions need to be computed. In addition, since the M2L matrices are practically approximated by low-rank ones \cite{FONG20098712} of the form $C[t,s] \approx \bo{L}[t,s]^*\bo{R}[t,s]$ with $\bo{L}[t,s],\bo{R}[t,s]\in \bb{R}^{k\times L^3}$, $k<<L^3$, we can exploit this factorization to further accelerate the M2L evaluation in Eq. \eqref{eq;mutualM2L}. This takes the following form:
\begin{equation}
\begin{aligned}
    \cc{M}_s^TC[s,t]\cc{M}_t &\approx \cc{M}_s^T\bo{L}[t,s]^*\bo{R}[t,s]\cc{M}_t\\
    &= \underbrace{\underbrace{\left(\bo{L}[t,s]\cc{M}_s\right)}_{\cc{O}(kL^3)}\underbrace{\left(\bo{R}[t,s]\cc{M}_t\right)}_{\cc{O}(kL^3)}.}_{\cc{O}(kL^3)}
\end{aligned}
\end{equation}
Different techniques can be used to obtain such low-rank approximations. We compared two of them: partial Adaptive Cross Approximation \cite{AIMI2022351} (pACA) and Singular Value Decomposition (SVD). pACA benefits from a linear algorithmic complexity while SVD is cubic cost. However, the sizes of the M2L matrices are relatively small and thanks to the use of symmetry \cite{chollet:hal-03563005}, it is possible to only compute a few dozens of such approximations per level (and only once during the precomputation step, not at runtime). Hence, because numerical ranks obtained using pACA being higher than ranks obtained through SVD (that are optimal), we measured a slightly higher cost when evaluating the FMM using pACA approximations. For these reasons, we decided to rely on SVDs for performance tests.

Similarly, because our kernels $K$ are transnationally invariant, we have
\begin{equation*}
    \kk{D}_\bo{x}\kk{D}_{\bo{y}}K(\bo{x},\bo{y}) + \kk{D}_\bo{y}\kk{D}_{\bo{x}}K(\bo{y},\bo{x}) = 2\kk{D}_\bo{x}\kk{D}_{\bo{y}}K(\bo{x},\bo{y}),
\end{equation*}
meaning that only half of the direct interactions have to be computed.

\subsection{Global ANKH-FMM algorithm}
For the sake of clarity, we summarize all the ANKH-FMM steps in Alg. \ref{alg:ankhfmm}. More precisely, ANKH-FMM first step consists in computing the modified charges in each leaf cell, once the (perfect) octree space decomposition is known. These modified charges are interpreted as terms of a IBFMM multipole expansion associated to interpolation nodes in these cells. This thus lead to a computation of other multipole expansions in the non-leaf tree levels. Notice that this step seems quite easy to parallelize in a shared-memory context due to the independence of the modified charges of two different leaves. The second ANKH-FMM step corresponds to the M2L evaluation involved in IBFMM, except that in our case mutual interactions give a simpler expression of it (following Sect. \ref{ss:mutual}). Since we are focused on energy computation and thanks to the adjoint form of the L2L/M2M and L2P/P2M operators, results of each mutual M2L can be added to a scalar variable. Thus, no L2L/L2P application is needed in practice.

An important point to mention at this presentation stage is that, even if the IBFMM part of ANKH-FMM also involves cells of the direct neighbor images of the simulation box $B$, there is no need of computing its multipole expansions. Indeed, there are equal to multipole expansions in the original box up a translation, which only appears in the M2L formula since all other formula are in practice taken respectively to cell centers \cite{chollet:tel-03203231}.

Then, the third step of ANKH-FMM is just a direct computation of the near field $\cc{N}_{real}$. To achieve it in practice, we used explicit derivations \cite{Smith2007PointMI} of the kernel $K$. Fourth step corresponds to the handling of periodicity (i.e. here the influence on $B$ of non-direct neighbor images of $B$) using the technique presented in Sect. \ref{ss:fmmperiodicity}. Finally, the two last steps are just direct computation of the self energy in Eq. \eqref{eq:ewald} and simple sum of all the computed quantities.
\begin{algorithm}[htbp]
\caption{ANKH-FMM}
\label{alg:ankhfmm}
\begin{algorithmic}
    \State // \textbf{\underline{Step 1:} get modified charges, multipole expansions}
    \For{each leaf cell $s\in \cc{L}(T)$}
        \State Compute $\cc{Q}_{s}$ using Alg. \ref{alg:compute_der_p_S}
        \State $\cc{M}_s = \cc{Q}_s$
    \EndFor
    \For{each non leaf level $e$ starting from the deepest one}
        \For{each cell $s$ at level $e$}
            \State $\cc{M}_s = \bo{0}$
            \For{each daughter $s'$ of $s$}
                \State $\cc{M}_s = \cc{M}_s + V[s,s']\cc{M}_{s'}$
            \EndFor
        \EndFor
    \EndFor
    \State // \textbf{\underline{Step 2:} get $\cc{F}_{real}$ using mutual M2L}
    \State $\cc{F}_{real} = \bo{0}$
    \For{each target cell $t$}
        \For{each source cell $s\in \Lambda(t)$}
            \If{interaction of $s$ and $t$ has not been yet computed}
                \State $\cc{F}_{real} = \cc{F}_{real}+2\cc{M}_t^TC[t,s]\cc{M}_s$
            \EndIf
        \EndFor
    \EndFor
    \State // \textbf{\underline{Step 3:} get $\cc{N}_{real}$ through mutual P2P}
    \State $\cc{N}_{real} = \bo{0}$
    \For{each leaf cell $t\in \cc{L}(\scr{T})$}
        \For{each source cell $s$ adjacent to $t$}
            \If{interaction of $s$ and $t$ has not been yet computed}
                \For{$\bo{x}\in X_{|t}$}
                    \For{$\bo{y}\in Y_{|s}$}
                        \State $\cc{N}_{real} = \cc{N}_{real}+2\kk{D}_\bo{x}\kk{D}_{\bo{y}}G(\bo{x},\bo{y})$
                    \EndFor
                \EndFor
            \EndIf
        \EndFor
    \EndFor
    \State // \textbf{\underline{Step 4:} get periodic influence}
    \State Let $\cc{M}_B$ be the root multipole expansion
    \State $\cc{F}_{per} = \cc{M}_B^T\chi^*\bb{F}^*D_B\bb{F}\chi \cc{M}_B$
    \State // \textbf{\underline{Step 5:} compute self energy}
    \State $\cc{E}_{self} = 0$
    \For{each $\bo{x}\in X$}
        \State $\cc{E}_{self} = \cc{E}_{self} + q_\bo{x}^2\frac{2\xi^2}{3}\mu_\bo{x}\cdot \mu_\bo{x} + \frac{8\xi^4}{5}\Theta_\bo{x}:\Theta_\bo{x}$
    \EndFor
    \State // \textbf{\underline{Step 6:} return energy $\cc{E}$}
    \State $\cc{E}=-\frac{\xi}{\sqrt{\pi}}\cc{E}_{self} + \cc{F}_{real} + \cc{N}_{real} + \cc{F}_{per}$
\end{algorithmic}
\end{algorithm}

\subsection{A first conclusion}
We presented in this section a new method, based on both the Ewald summation and the interpolation-based FMM, aiming at efficiently treating $N$-body problems arising in the energy computation context in molecular dynamics and in a way that should efficiently perform in parallel (mainly distributed memory) context. Hence, since FMMs benefit from highly scalable strategies \cite{10.5555/147877.148090,10.1145/2160718.2160740,10.1145/1654059.1654123}, this new ANKH-FMM approach may pave the way for a linear complexity and scalable family of methods suited for molecular large scale dynamics computations. We proved, under molecular dynamics hypotheses, that our method is actually $\cc{O}(N)$, even when exploiting periodicity with large amount of images, thanks to a new way of taking into account the influence of the far images. We introduced the mutual interactions in the context of M2L evaluations for this application case and we presented a numerical differentiation through polynomial interpolation in order to deal with differential operators arising when considering dipoles and quadrupoles (our methodology can be extended to any multipole order with low effort). Finally, we provided explicit algorithm detailing our method.

Among the direct perspective, the ideas behind ANKH-FMM can be used to compute forces, even if some optimizations have to be changed to this purpose (especially in the mutual M2L interactions). The parallel implementation of ANKH-FMM as well as real tests on supercomputers also counts among the important and relatively short-term perspectives.


\section{ANKH-FFT}
\label{s:ankhfft}
Despite the sequential and distributed efficiency of the interpolation-based FMM, limitations arise in the context of GPU computing. Indeed, if the near-field part is well suited to GPU, the M2L evaluation suffers from data loading in shared memory that mitigates the performance on this type of architecture.
However, the perfectness of the $2^d$-trees in molecular dynamics allows to derive efficient schemes, as those proposed in the literature \cite{gpudarve}. In this section, we introduce a new method, also based on interpolation, and designed to bypass the limitations arising in GPU context.

Starting from the interpolation of the particle multipoles at leaves on equispaced grids (instead of Chebyshev ones), one may construct a single level FMM as the evaluation of the matrix concatenating \textbf{all} the M2L matrices between non-adjacent or equal cells. This matrix can be proved to be a block-Toeplitz matrix (we will focus on the theory in a forthcoming mathematical article). Such a matrix can be embedded in a larger circulant matrix. Because any circulant matrix can be explicitly diagonalized in the Fourier domain, and because the product by a discrete Fourier basis can be efficiently processed through FFTs, one may compute this diagonalization in a linearithmic time. Doing so, we obtain a fast linearithmic scheme to deal with the far-field of N-body problems on quasi-uniform particle distributions (this assumption allowing to obtain perfect trees). Of course, the near-field is treated as in the FMM, i.e. by performing local direct computation between non-separated cells.

However, in our case, since we want to perform interpolation on Chebyshev polynomials in order to mitigate the precision loss of the multipole operator evaluation without much numerical stability issue, we start by performing such a Chebyshev interpolation at leaves, directly followed by a re-interpolation of the data now defined at Chebyshev nodes to a local equispaced grid, using the same grid (up to translation) for every leaf (since they all have the same radius).

Once in the Fourier domain, the energy computation problem is nothing more than a scalar product between two vectors: the diagonal of the Fourier diagonalization of the circulant matrix and the modified charges vector Fourier transform squared modulus after zero-padding (to have the same size than the circulant matrix). Since the circulant matrix has $\cc{O}(N)$ rows and columns ($\cc{O}(N)$ leaves with $\cc{O}(1)$ interpolation nodes in each), this scalar product costs $\cc{O}(N)$ floating point operations (flops) while the FFT costs $\cc{O}(NlogN)$ flops. Since the interpolation phase is analogue to the application of all P2M operators in ANKH-FMM, this step can also be performed in $\cc{O}(N)$ flops. Hence, we end up with a linearithmic method. This new method is named ANKH-FFT.

ANKH-FFT is theoretically more costly than ANKH-FMM, but ANKH-FFT benefits from highly optimized FFT implementations, also available on GPU, and the interpolation phase is a highly parallel one. Hence, ANKH-FFT seems to be a better candidate than ANKH-FMM in a GPU context. We should still compare the performance of these two methods, which is done in the following. These performance are comparable, and ANKH-FFT behaves as a linear method in the tested range (i.e. from 1000 to 100 Million atoms). This leads us to the conclusion about the use of this new approach: ANKH-FFT appears as an entirely new and efficient approach to deal with local $N$-body problems with an algorithmic structure well-suited for GPU programming. This should thus be an excellent alternative to local interpolation-based FMM used in one core of the distributed memory framework of the parallel interpolation-based FMM.

We may here insist on the fact that ANKH-FFT requires a strong mathematical background to be entirely presented, hence we postpone this presentation to a forthcoming mathematical paper, including proofs. ANKH-FFT is based on a matrix factorization presented in Sect. \ref{subsubsect:matrixfact}. We then show how this factorization leads to a diagonalization in Sect. \ref{ss::ankh_fft_diag} and we present the fast periodic condition handling in Sect. \ref{ss::effperiofft}. The important role of mutual interactions in ANKH-FFT is presented in Sect. \ref{ss::ankhfftmutual}. Finally, we summarize the main algorithmic in Sect. \ref{ss::ankhdfft_globalalgo}, followed in Sect. \ref{ss::ankhfft_complexity} by a complexity analysis. The numerical comparison between ANKH-FFT and ANKH-FMM is provided in Sect. \ref{ss::comparisonfftfmm}.

\subsection{Matrix factorization}
\label{subsubsect:matrixfact}
\begin{figure}
    \centering
    \includegraphics[width=\linewidth]{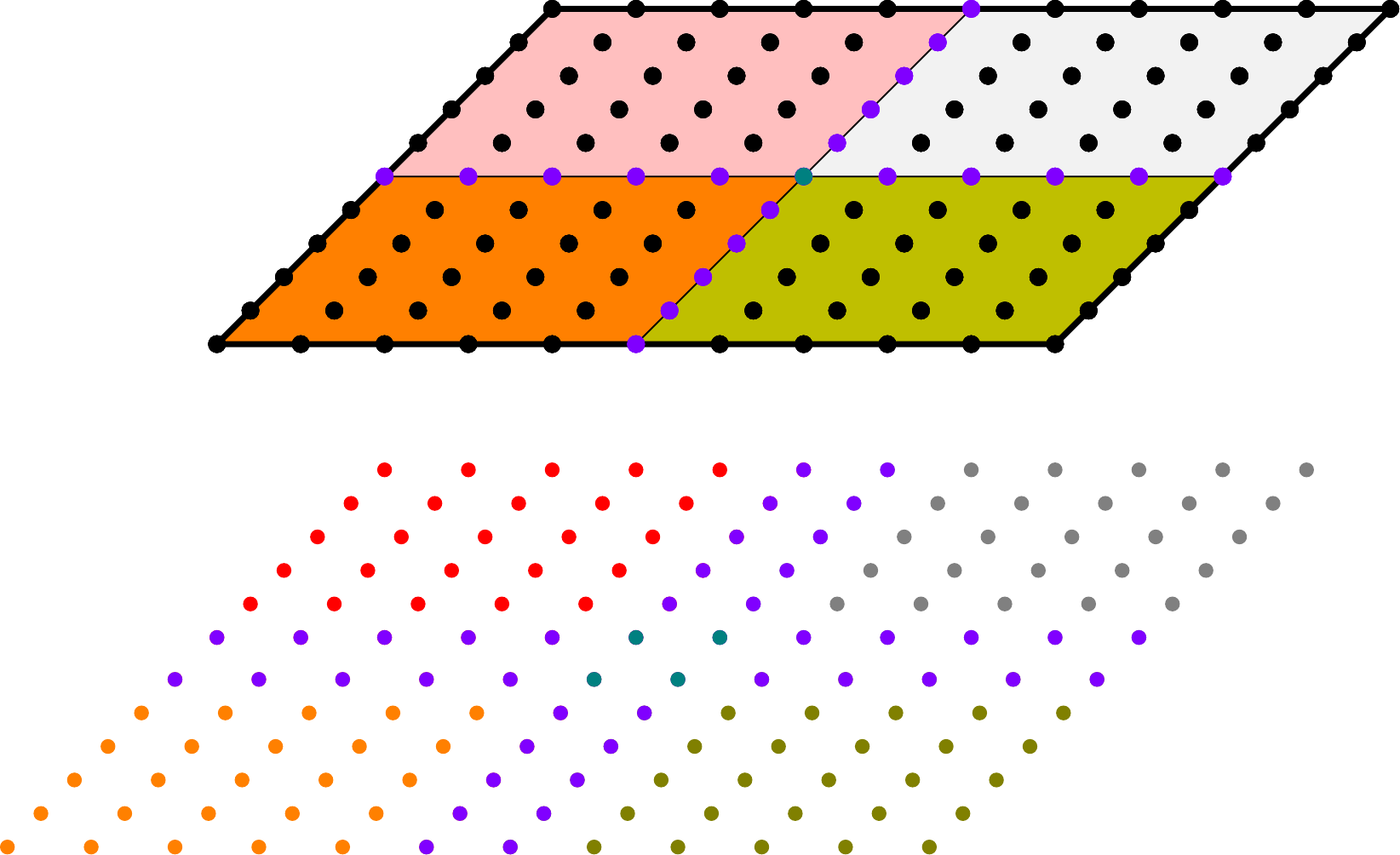}
    \caption{Schematic representation of 2D interpolation on tensorized equispaced grids (as in Eq. \ref{eq::tensorized_equispaced_grid}) with $L=5$ on four different leaf cells (in red, orange, green and grey). Interpolation nodes belonging on a unique cell are represented in black on the top figure, nodes belonging in two cells are represented in purple and the node belonging in the four cells is represented in cyan. In the bottom, we represented the actual grid used in AKNH-FFT, with a slight node translation in order to exhibit duplications: purple nodes are duplicated once and cyan one is duplicated four times. Assuming that we only have these four leaf cells (i.e. a tree depth of only one), this cyan node is indexed (following the multi-indexing of Sect. \ref{ss::ankh_fft_diag}) as $(0,4,0,0)$ in the red cell, as $(1,0,0,0)$ in the grey cell, as 
$(0,4,0,4)$ in the orange cell and as $(1,0,1,4)$ in the green cell.\label{fig:allequispacedgrids}}
\end{figure}

We start by considering all the leaf cells $\cc{L}(\scr{T})$ of $\scr{T}$ (that are all at the same tree level since $\scr{T}$ is perfect). Let us consider the \textit{extended interaction list} $\Lambda^+$ of any leaf cell $t$ defined as the set of leaf cells with a strictly positive distance with $t$:
\begin{equation*}
    \Lambda^+(t) := \{s\in \cc{L}(\scr{T})\hspace{0.1cm}|\hspace{0.1cm}dist(t,s) > 0\}.
\end{equation*}
We then decompose $\tilde{\cc{E}}_{real}$ of Eq. \ref{eq_real_nbody} as
\begin{equation}
\begin{aligned}
    \cc{\tilde{E}}_{real} =& \underbrace{\sum_{t\in \cc{I}(\scr{T})}\sum_{\substack{\bo{x}\in t,\bo{y}\in s\\s\in \Lambda^+(t)}}\kk{D}_\bo{x}\kk{D}_\bo{y}H(\bo{x},\bo{y})}_{=:\cc{F}}  \\&+ \underbrace{\sum_{t\in \cc{I}(\scr{T})}\sum_{\substack{\bo{x}\in t,\bo{y}\in s\\s\notin \Lambda^+(t)}}\kk{D}_\bo{x}\kk{D}_\bo{y}H(\bo{x},\bo{y})}_{=:\cc{N}}.
\end{aligned}
\end{equation}
As in our ANKH-FMM algorithm, we propose to compute $\cc{N}$ directly and to approximate $\cc{F}$ by means of interpolation. Indeed, following Sect. \ref{ss:mutual}, we can write
\begin{equation}
    \cc{F} \approx \sum_{t\in \cc{I}(\scr{T})}\sum_{s\in \Lambda^+(t)}\cc{M}_t^TC[t,s]\cc{M}_s.
\end{equation}
Now, let $\bo{E}_{c}$, for any $c\in \cc{L}(\scr{T})$ be a reinterpolation matrix (see the M2M definition in Eq. \ref{eq:M2M}) from the (Chebyshev) interpolation nodes of $c$ into a equispaced grid $\cc{J}_L$ defined by
\begin{equation}
\label{eq::tensorized_equispaced_grid}
\begin{aligned}
    \cc{J}_L := \Bigg\{&r_B2^{-(1+E)}\left(\begin{bmatrix}
      -1\\-1\\-1
    \end{bmatrix}+\frac{2\bo{k}}{L-1}\right)\hspace{0.1cm}\Bigg|\hspace{0.1cm}\\&\bo{k}\in [\![0,L-1]\!]^3\Bigg\},
    \end{aligned}
\end{equation}
i.e. $\bo{E}_c$ interpolates the modified charges associated to Chebyshev nodes in a cell $c$ to new modified charges associated to the nodes of a equispaced grid in $c$. This interpolation grid $\cc{J}_L$ has to be centered in each leaf cell center. Said differently, the same equispaced interpolation grid $\cc{J}_L$ is used in each leaf cell, up to a translation (see Fig. \ref{fig:allequispacedgrids}).

Each M2L matrix $C[t,s]$ can be written as
\begin{equation}
\label{eq::fromchebtoequispaced}
    C[t,s] \approx \bo{E}^T_tG[t,s]\bo{E}_s
\end{equation}
where $\left(G[t,s]\right)_{k,l} := H(ctr(t)+\tilde{\bo{x}}_k,ctr(s)+\tilde{\bo{x}}_l)$, $ctr(c)$ denoting the center of the cell $c$ and $\tilde{\bo{x}}_k$ referring to the $k^{th}$ elements of $\cc{J}_L$. In other terms, Eq. \ref{eq::fromchebtoequispaced} simply switches between Chebyshev nodes in a cell and equispaced nodes.

\begin{remark}
\label{rem::compression}
   In practice, the order of the Chebyshev rules are chosen so that the derivation of the associated polynomials results in the targeted error. Nevertheless, after modified charges computation, the equispaced interpolation does not have to consider the same interpolation order. This last order can be chosen without consideration about the derivatives, which practically implies lower orders. Hence, another compression level may be induced by this reinterpolation, provided that the interpolation order on equispaced grids is chosen to be smaller than the one on Chebyshev grids.
\end{remark}

Now, we aim at efficiently computing the interactions between a cell $t$ and all the cells in its extended interaction list. Let us extend the definition of $G[t,s]$ in order to cover the blanks induced by pairs of leaf cells $(t,s)$ such that $s\notin\Lambda^+(\scr{T})$, introducing
\begin{equation}
\label{eq::gtilde}
    \tilde{G}[t,s] = \begin{cases}
    G[t,s]\textit{ if }s\in\Lambda^+(t),\\
    \bo{0}\textit{ otherwise}.
    \end{cases}
\end{equation}
We thus obtain 
\begin{equation}
\begin{aligned}
    \cc{F} &\approx \sum_{t\in \cc{L}(\scr{T})}\sum_{s\in \cc{L}(\scr{T})}\cc{M}_t^T\bo{E}^T_t\tilde{G}[t,s]\bo{E}_s\cc{M}_s
    \end{aligned}
\end{equation}
that corresponds to the product of a matrix $\bb{A}$ formed by the concatenation of the blocks $\tilde{G}[t,s]$ with a left and a right vector formed by the concatenation of $\bo{E}_c\cc{M}_c$'s, $c\in \cc{L}(\scr{T})$. The blocks are ordered in 3D (in both rows and column) following a lexicographic order on their center difference:
\begin{equation*}
\begin{aligned}
    c < c' \Leftrightarrow \Big[r_0 < 0 &\hspace{0.1cm}\lor\hspace{0.1cm} \left(r_0 = 0 \land r_1 < 0\right)\\&\hspace{0.1cm}\lor\hspace{0.1cm} \left(r_0 = 0 \land r_1 = 0\land r_2 < 0\right)\Big]
\end{aligned}
\end{equation*}
for any $c,c'\in \cc{L}(\scr{T})$ and with $(r_0,r_1,r_2) := ctr(c)-ctr(c')$. The result of this concatenation is a vector $\bo{v}\in \bb{R}^{8^EL^3}$, where $8^E$ corresponds to the number of leaf cells in a perfect octree with depth $E$. Using these notations, we obtain the following approximation:
\begin{equation}
\label{eq::fasvav}
    \cc{F} \approx \bo{v}^T\bb{A}\bo{v}.
\end{equation}
We are thus now interested by an efficient computation of the product by $\bb{A}$. This is the purpose of Sect. \ref{ss::ankh_fft_diag}.

\begin{remark}
In $\bb{A}$, a given target cell $t$ may interact with various source cells $s$ that spatially share interpolation nodes. However, each equispaced interpolation grid and each multipole expansion $\cc{M}_s$ associated to it has to be considered independently because, depending on $t$, $\cc{M}_s$ may be masked or used, thanks to possible products by zeros in Eq. \ref{eq::gtilde}. This means that we do not see the concatenation of equispaced grids as a global equispaced grid, but rather as a collection of such small grids with nodes duplication on the cell's edges. This consideration is illustrated on Fig. \ref{fig:allequispacedgrids}
\end{remark}

\subsection{Diagonalization}
\label{ss::ankh_fft_diag}
In this section, we provide an explicit diagonalization of an embedding of $\bb{A}$ in a Fourier basis. Let $\bo{k},\bo{l}$ be multi-indices in
\begin{equation*}
\begin{aligned}
    \left(\underbrace{[\![0,n(E)]\!]}_{\substack{\text{cell on}\\\text{ x-axis}}}\times\underbrace{[\![0,L-1]\!]}_{\substack{\text{interpolation}\\\text{node on}\\\text{ x-axis}}}\right)&\times
    \left(\underbrace{[\![0,n(E)]\!]}_{\substack{\text{cell on}\\\text{ y-axis}}}\times\underbrace{[\![0,L-1]\!]}_{\substack{\text{interpolation}\\\text{node on}\\\text{ y-axis}}}\right)\\&\times
    \left(\underbrace{[\![0,n(E)]\!]}_{\substack{\text{cell on}\\\text{ z-axis}}}\times\underbrace{[\![0,L-1]\!]}_{\substack{\text{interpolation}\\\text{node on}\\\text{ z-axis}}}\right),
\end{aligned}
\end{equation*}
where $n(E) := 2^E-1$ (see Fig. \ref{fig:allequispacedgrids} for graphical 2D representation). Let $\bo{u}=(\tilde{u}_0,\tilde{u}_1,\tilde{u}_2)\in\{0,1\}^6$ be a boolean vector and let $\tilde{\bo{k}}=(\tilde{k}_0,\tilde{k}_1\tilde{k}_2),\tilde{\bo{l}}=(\tilde{l}_0,\tilde{l}_1,\tilde{l}_2)$ be defined as
\begin{equation}
    \tilde{k}_p = \begin{cases}k_p + u_p\textit{ if }p\textit{ is even and }k_p\neq 2^E-1,\\k_p + u_p\textit{ if }p\textit{ is odd and }k_p\neq L-1,\\k_p\textit{ otherwise}.\end{cases}
\end{equation}

In our study case, $\bb{A}$ verifies a strong property which is given in Thm. \ref{thm:main_theo}.

\begin{theorem}
\label{thm:main_theo}
If $H$ is a radial kernel, then $\bb{A}_{\bo{k},\bo{l}} = \bb{A}_{\tilde{\bo{k}},\tilde{\bo{l}}}$, where $\bb{A}_{\bo{k},\bo{l}}$ refers to the element of $\bb{A}$ at row $k_0(2^EL^3)^2+k_12^EL^3+k_2$ and column $l_0(2^EL^3)^2+l_12^EL^3+l_2$.
\end{theorem}
Actually, Thm. \ref{thm:main_theo} exhibits a Toeplitz structure for $\bb{A}$ in six dimensions (that are the dimensions of the indexing at the begining of this section). The direct consequence of Thm \ref{thm:main_theo} is summarized in Cor. \ref{cor:main_cor}.
\begin{corollary}
\label{cor:main_cor}
Let $\chi \in \{0,1\}^{ \left(\left((2^{E+1}-1)\times (2L-1)\right)^3\right)\times \left(\left(2^E\times L\right)^3\right)}$ be defined as
\begin{equation}
    \chi_{\bo{k},\bo{l}} = \begin{cases}1\textit{ if }\bo{k}=\bo{l},\\
    0\textit{ otherwise}.\end{cases}
\end{equation}
for any $\bo{k}\in \left([\![0,2^E-1]\!]\times[\![0,L-1]\!]\right)^3,\bo{l}\in \left([\![0,2^{E+1}-2]\!]\times[\![0,2L-2]\!]\right)^3$ and let $\bb{F}_6$ be the discrete Fourier transform matrix of dimensions $(2^E,L,2^E,L,2^E,L)$. We have
\begin{equation}
\label{eq::diag_of_A}
    \bb{A} = \chi^*\bb{F}_6^*\bo{D}(\bb{A})\bb{F}_6\chi,
\end{equation}
where $\bo{D}(\bb{A})$ is a diagonal matrix.
\end{corollary}
The main idea is that $\bb{A}$ is a block-Toeplitz matrix according to Thm. \ref{thm:main_theo}. Such a matrix can be embedded into a block-circulant matrix one. Any block-circulant matrix can be diagonalized in a Fourier basis (here $\bb{F}_6$), giving the factorization of Cor. \ref{cor:main_cor}. Here, this block-circulant matrix is never built because the diagonal of $\bo{D}(\bb{A})$ can be computed through a single application of $\bb{F}_6$ to a vector of correctly ordered elements of the first row and column of $\bb{A}$. Similar tools were used in Sect. \ref{ss:fmmperiodicity} in a much simpler context.

The expression \eqref{eq::diag_of_A} allows to derive a fast scheme for the product by $\cc{F}$ in Eq. \eqref{eq::fasvav} since thanks to Thm. \ref{thm:main_theo}:
\begin{equation}
\label{eq::fastotfft}
    \cc{F} \approx \bo{v}^T\chi^*\bb{F}_6^*\bo{D}(\bb{A})\bb{F}_6\chi\bo{v},
\end{equation}
in which products by $\chi^*, \chi$ can be evaluated with linear complexity (because of their sparse structure) and products by $\bb{F}_6^*,\bb{F}_6$ can be processed with linearithmic complexity through FFTs (regarding the size of $\bo{v}$).

We refer to Alg. \ref{alg:gen_first_col} for the explicit computation of the diagonal of $\bo{D}(\bb{A})$. In this algorithm, we used the notations $\bo{I} = (I_0,...,I_5)^T$, $\bo{J} = (J_0,...,J_5)^T$, $\bo{R} = (R_0,...,R_5)^T$. Mainly, this algorithm generates a vector composed of the (correctly ordered for the circulant embedding) elements of the first row and column of $\bb{A}$, followed by the discrete Fourier transform of this vector. This results in the diagonal elements of the diagonalization of $\bb{A}$ in the Fourier domain, i.e. $\bo{D}(\bb{A})$.

\begin{algorithm}[htbp]
\caption{Get ANKH-FFT diagonal matrix}
\label{alg:gen_first_col}
\begin{algorithmic}
\State $u = \frac{(2^{E}-1)r_B}{2^{E}}$
\State $v = \frac{r_B}{2^E}$
\State $P_{cell} = 2^{E+1}-1$
\State $P_{ntrp} = 2L-1$
\State $\bo{Z} := \begin{bmatrix}
  \frac{2u}{2^E-1} & \frac{2v}{L-1} & 0 & 0 & 0 & 0\\
  0 & 0 & \frac{2u}{2^E-1} & \frac{2v}{L-1} & 0 & 0 \\
  0 & 0 & 0 & 0 & \frac{2u}{2^E-1} & \frac{2v}{L-1}
\end{bmatrix}$
\For{$\bo{R}\in \left([\![0,P_{cell}-1]\!]\times [\![0,P_{ntrp}-1]\!]\right)^3$}
    \For{$k=0,1,2$}
        \If{$R_{2k} < 2^E$}
            \State $I_{2k} = R_{2k}$ and $J_{2k} = 0$
        \Else 
            \State $I_{2k} = 0$ and $J_{2k} = P_{cell}-R_{2k}$
        \EndIf
        \If{$R_{2k+1} < L$}
            \State $I_{2k+1} = R_{2k}$ and $J_{2k+1} = 0$
        \Else 
            \State $I_{2k+1} = 0$ and $J_{2k+1} = P_{cell}-R_{2k}$
        \EndIf
        \State $\bo{x} = \bo{Z}\bo{I}$
        \State $\bo{y} = \bo{Z}\bo{J}$
    \EndFor
    \State $t = 0$
    \State $n = d$
    \For{$k=0,1,2$}
        \State $t = t+|I_{2k} - J_{2k}|$
        \If{$(I_{2k} - J_{2k}) \neq 0$}
            \State $n = n-1$
        \EndIf
    \EndFor
    \State $i=\sum_{k=0}^{2}(R_{2k}P_{ntrp}+R_{2k+1})P_{cell}^{2-k}P_{ntrp}^{2-k}$
    \If{$t > n$}
        \State $\bo{v}_{i} = \sum_{\bo{t}\in 2r_B Z_p}H\left(\bo{x},\bo{y}+\bo{t}\right)$
    \Else
        \State $\bo{v}_{i} = 0$
    \EndIf
\EndFor
\State $diag\left(\bo{D}(\bb{A})\right) = \bb{F}_6\bo{v}$
\end{algorithmic}
\end{algorithm}

\subsection{Efficient handling of the periodicity}
\label{ss::effperiofft}
When large number of periodic images are considered, the computation of $\sum_{\bo{t}\in 2r_B Z_p}H\left(\bo{x},\bo{y}+\bo{t}\right)$ involved in Alg. \ref{alg:gen_first_col} becomes too coslty to be used. Fortunately, this step can be accelerated through another level of interpolation. Exploiting translational invariance of $H$, one gets:
\begin{equation}
    \begin{aligned}
    \sum_{\bo{t}\in 2r_B Z_p}H(\bo{x},\bo{y}+\bo{t}) &= \sum_{\bo{t}\in 2r_B Z_p}H(\underbrace{\bo{x}-\bo{y}}_{=:\bo{z}},\bo{t})
    \end{aligned}
\end{equation}
that can be separated into two sums
\begin{equation}
\label{eq::splitper}
    \underbrace{\left(\sum_{\bo{t}\in 2r_B Z_2}H(\bo{z},\bo{t})\right)}_{=:\cc{N}_{per}(\bo{z})} + \underbrace{\left(\sum_{\bo{t}\in 2r_B Z_p\backslash Z_2}H(\bo{z},\bo{t})\right)}_{\cc{F}_{per}(\bo{z})}.
\end{equation}
Since both $\bo{x},\bo{y}\in B$ lie inside the (centered at zero) simulation box of radius $r_B$, $\bo{z}$ belongs inside a box $\tilde{B}$ of radius $2r_B$. Each $\bo{t}\in 2r_B Z_p\backslash Z_2$ has distance at least $2r_B$ from $\tilde{B}$. hence, each such $\bo{t}$ is well-separated from $\tilde{B}$, using a treecode-like criterion \cite{Dehnen_2002}. Hence, performing an interpolation over a equispaced grid on $\tilde{B}$, we obtain
\begin{equation}
\label{eq::interpfper}
    \cc{F}_{per}(\bo{z}) \approx \sum_{k}U_k[\tilde{B}](\bo{z})\cc{F}_{per}(\bo{z_k}),
\end{equation}
$U_k$ being the Lagrange polynomial associated to the $k^{th}$ node of the equispaced interpolation grid over $\tilde{B}$.

Combining Eq. \eqref{eq::splitper} and Eq. \eqref{eq::interpfper}, we get
\begin{equation}
\label{eq::nbodyforperio}
\begin{aligned}
    \cc{F}_{per}(\bo{z}) &\approx \sum_kU_k[\tilde{B}](\bo{z})\underbrace{\left(\sum_{\bo{t}\in 2r_B Z_p\backslash Z_2}H(\bo{z}_k,\bo{t})\right)}_{=:h_k}
\end{aligned}
\end{equation}
in which the sum into parentheses corresponds to a $N$-body problem on two point clouds: the Caresian grid over $\tilde{B}$ (with $L^3$ nodes) and the set $2r_B Z_p\backslash Z_2$ (with all charges equal to $1$). This $N$-body problem can be solved efficiently using a treecode or a IBFMM. Then, once all $h_k$'s are known, any $\cc{F}_{per}(\bo{z})$ can be retrieve in $\cc{O}(L^3) = \cc{O}(1)$. In addition, these $h_k$'s can be computing once during the precomputation step since they do not depend on the modified charges but only on the kernel $H$.

The term $\cc{N}_{per}(\bo{z})$ involves only a sum over $5^3 = \cc{O}(1)$ elements, meaning that using this interpolation, one may deduce the result of Eq. \eqref{eq::splitper} in $\cc{O}(1)$ flops.

\subsection{Mutual interactions}
\label{ss::ankhfftmutual}
The evaluation of Eq. \eqref{eq::fastotfft} involves two FFTs if computed naively. Since they are the most consuming part of the evaluation time, we may wonder how to mitigate it. This can be done by exploiting once again, as in Sect. \ref{ss:mutual}, the mutual interactions. Our idea is based on the following result:
\begin{corollary}
\label{cor::realspectrum}
    If $H$ has real values, then $\bo{D}(\bb{A})\in \bb{R}^{[(2^{E+1}-1)(2L-1)]^3\times [(2^{E+1}-1)(2L-1)]^3}$.
\end{corollary}
This can be proved from Thm. \ref{thm:main_theo} noticing that that radial property of $H$ implies symmetry of $\bb{A}$. Hence, the spectrum of $\bb{A}$ is real.

As a direct consequence, following Eq. \eqref{eq::fastotfft}, we have
\begin{equation}
\begin{aligned}
    \bo{v}^T\chi^*\bb{F}_6^*\bo{D}(\bb{A})\bb{F}_6\chi\bo{v} &= \left(\bb{F}_6\chi\bo{v}\right)^*\bo{D}(\bb{A})\underbrace{\bb{F}_6\chi\bo{v}}_{=:\bo{w}},
\end{aligned}
\end{equation}
with $\bo{w}=(w_0,...,w_{[(2^{E+1}-1)(2L-1)]^3-1})^T$, and which reformulates as
\begin{equation}
\label{eq::realdiag}
\begin{aligned}
    \bo{w}^*\bo{D}(\bb{A})\bo{w}&=
    \sum_{k=0}^{[(2^{E+1}-1)(2L-1)]^3-1}\overline{w_k}\bo{D}(\bb{A})_{k,k}w_k\\
    &= \sum_{k=0}\bo{D}(\bb{A})_{k,k}|w_k|^2.
\end{aligned}
\end{equation}
Hence, only $\bo{w}$ has to be computed, meaning that only a single FFT needs to be performed during evaluation of Eq. \eqref{eq::fastotfft}. As another important consequence of Cor. \ref{cor::realspectrum}, once this FFT has been performed, one can only keeps the real part of the output (the imaginary one being equal to 0). In our implementation, we even rely on real-to-complex FFTs from FFTW3 \cite{FFTW05} since the modified charges (as well as the entries of $\bo{w}$) are real provided that the charges, dipoles and quadrupoles also are.

\subsection{Global algorithmic}
\label{ss::ankhdfft_globalalgo}
It is now possible to provide the overall algorithm for the method derived from our ANKH methodology and exploiting the diagonalization of the interpolated Ewald summation over local equispaced grids. Namely, this corresponds to ANKH-FFT and Alg. \ref{alg:ankhfft} summarizes all its steps.

ANKH-FMM and ANKH-FFT shares strong similarities: both use the same computation of modified charges, both compute the near field part ($\cc{N}_{real}$) using direct interactions and both compute the the self energy in the same way. However, the computation of the far field part ($\cc{F}_{real}$) is strongly simplified in ANKH-FFT.

Among the differences, we decided to take into account the influence of the periodicity in the diagonalization process of ANKH-FFT while we opted for an external handling in ANKH-FMM, meaning that the evaluation of ANKH-FFT does not involve any periodic influence (already contained in $\bo{D}(\bb{A})$). The reason behind this choice is that in ANKH-FFT, the far field computation is done in the Fourier domain directly, in which the particle cannot be interpreted as easily as in the cartesian domain, while in ANKH-FMM, the multilevel structure provide a way of measuring the influence of the far images efficiently \textit{on-the-fly}. This points out the main difference between the two approaches: ANKH-FFT is not multilevel but is built on a global operation (FFT).

\begin{algorithm}[htbp]
\caption{ANKH-FFT}
\label{alg:ankhfft}
\begin{algorithmic}
    \State // \textbf{\underline{Step 1:} get modified cartesian multipole expansions}
    \For{each leaf cell $s\in \cc{L}(T)$}
        \State Compute $\cc{Q}_{s}$ using Alg. \ref{alg:compute_der_p_S}
        \State $\cc{M}_s = \bo{E}_s\cc{Q}_s$
    \EndFor
    \State // \textbf{\underline{Step 2:} get $\cc{F}_{real}$ using FFT}
    \State Let $\bo{v}$ be defined as in Eq. \eqref{eq::fastotfft}
    \State $\bo{w} = \bb{F}_6\chi\bo{v}$
    \State $\cc{F}_{real} = \bo{w}^T\bo{D}(\bb{A})\bo{w}$ using Eq. \eqref{eq::realdiag}
    \State // \textbf{\underline{Step 3:} get $\cc{N}_{real}$ through mutual P2P}
    \State $\cc{N}_{real} = \bo{0}$
    \For{each leaf cell $t\in \cc{L}(\scr{T})$}
        \For{each source cell $s$ adjacent to $t$}
            \If{interaction of $s$ and $t$ has not been yet computed}
                \For{$\bo{x}\in X_{|t}$}
                    \For{$\bo{y}\in Y_{|s}$}
                        \State $\cc{N}_{real} = \cc{N}_{real}+2\kk{D}_\bo{x}\kk{D}_{\bo{y}}G(\bo{x},\bo{y})$
                    \EndFor
                \EndFor
            \EndIf
        \EndFor
    \EndFor
    \State // \textbf{\underline{Step 4:} compute self energy}
    \State $\cc{E}_{self} = 0$
    \For{each $\bo{x}\in X$}
        \State $\cc{E}_{self} = \cc{E}_{self} + q_\bo{x}^2\frac{2\xi^2}{3}\mu_\bo{x}\cdot \mu_\bo{x} + \frac{8\xi^4}{5}\Theta_\bo{x}:\Theta_\bo{x}$
    \EndFor
    \State // \textbf{\underline{Step 5:} return energy $\cc{E}$}
    \State $\cc{E}=-\frac{\xi}{\sqrt{\pi}}\cc{E}_{self} + \cc{F}_{real} + \cc{N}_{real}$
\end{algorithmic}
\end{algorithm}

\subsection{Complexity}
\label{ss::ankhfft_complexity}
There are two complexities to count in ANKH-FFT: the precompuation complexity as well as the evaluation one. The precomputation step involves the solving of $N$-body problem of Eq. \eqref{eq::nbodyforperio}, which may be done in a linear time with respect to the number of far images, i.e. in $\cc{O}(p^3)$. In practice, $p=\cc{O}(1)$ is small. Then, the diagonal of $\bb{A}$ is computed using Alg. \ref{alg:gen_first_col}. It is easy to see from this algorithm that, thanks to Sect. \ref{ss::effperiofft}, its complexity is linearithmic with respect to the size of $\bo{v}$.

When considering quasi-uniform particle distributions, the total number of leaf cells can be considered as $\cc{O}(N)$. In this case, the size of $\bo{v}$ is $\cc{O}(NL^3)$ with a constant $L$ corresponding to the number of interpolation node in one direction of the equispaced grid, which is a small constant. Hence, the size of $\bo{v}$ is $\cc{O}(N)$ and the FFTs perform in $\cc{O}(NlogN)$. 

The precomputation cost of the modified charges being linear with respect to the number of particles the simualtion box, the overall precomputation complexity is
\begin{equation}
    \cc{O}(N)+\cc{O}(N)+\cc{O}(NlogN) = \cc{O}(NlogN).
\end{equation}

Since the product by a diagonal matrix is also $\cc{O}(N)$, the overall complexity of the evaluation of $\cc{F}$ complexity though Eq. \eqref{eq::fastotfft} is $\cc{O}(NlogN)$ (because of the FFT).

\subsection{Comparison between ANKH-FMM and ANKH-FFT and partial conclusions}
\label{ss::comparisonfftfmm}
\begin{figure}
    \centering
    \includegraphics[width=\linewidth]{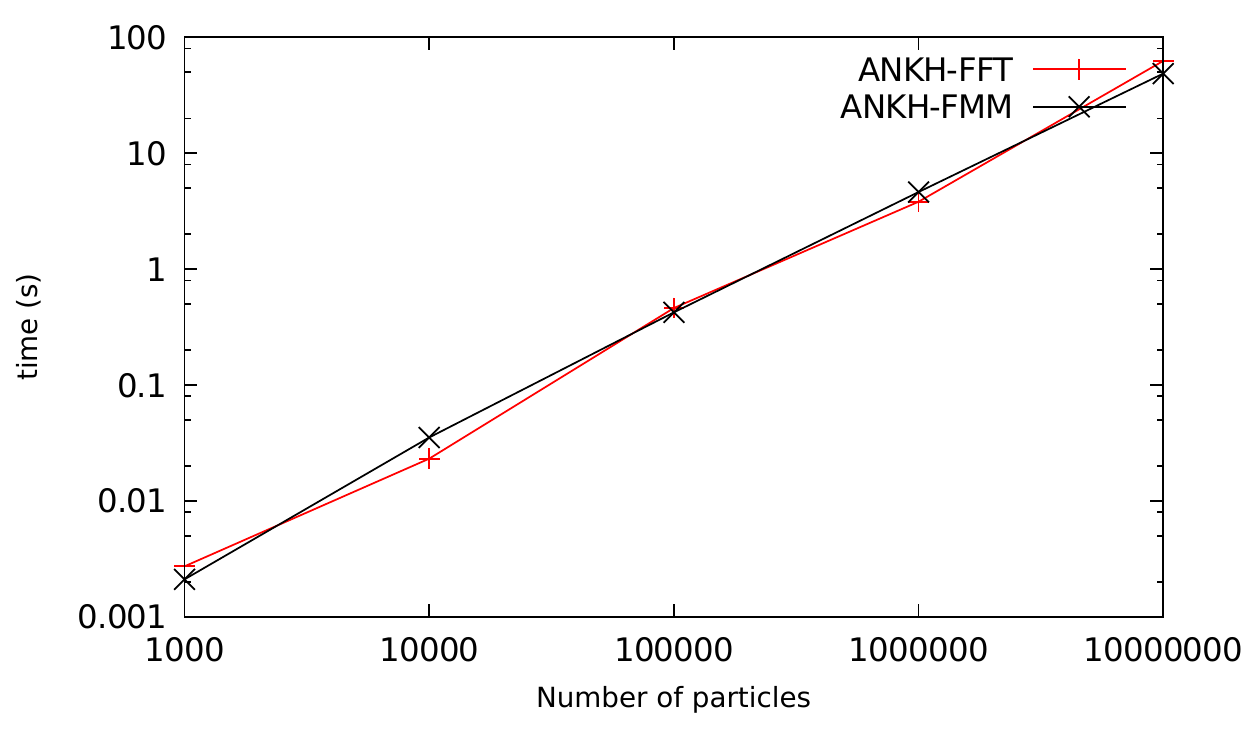}
    \caption{Timing comparison of the two ANKH variants with respect to the number of particles (atoms) in the simulation box.\label{fig:perfs_fmm_fft}}
\end{figure}

This section is dedicated to the comparison between ANKH-FMM and ANKH-FFT. As we already discussed, the theoretical complexities (respectively linear and linearithmic) are not sufficient to conclude in the practical efficiency. One of the main reasons behind that is the arithmetic intensity of the \textit{state-of-the-art} FFT implementations \cite{FFTW05} that should be higher than our IBFMM one. Hence, we compared running times for the evaluation of both methods (i.e. excluding precomputations) on the same particle distributions, using charges only (since the near field part is the same between the two methods) and no periodic condition (since its handling has a very small impact on the timings and that we want to compare the main core of these methods). Timings are provided in Fig. \ref{fig:perfs_fmm_fft}. To do these tests, we run our (sequential) codes on a single Intel(R) Xeon(R) Gold 6152 CPU.

Performance of the two methods are comparable in the (large) tested range of systems, i.e. from 1000 atoms to 10M atoms: ANKH-FFT performs slightly better on the case with 10000 atoms, but not in a significant way. Since ANKH-FMM follows the theoretical estimates of linear complexity, ANKH-FFT also behaves as a linear method in practice.

We may still insist on the fact that the precomputations involved in ANKH-FFT are more costly than the ones of ANKH-FMM (but in the same order of magnitude than the evaluation step). These precomputation timings are however quite relative to the application: in the MD context, the same simulation box is used for a massive amount of evaluations, hence these precomputations are performed only once for many evaluations, which mitigates their costs.

To summarize, we list the main difference, advantages and drawbacks of the two ANKH approaches:

\begin{itemize}
    \item ANKH-FMM is a truly linear method while ANKH-FFT only is a linearithmic one. Both achieve fast computations with comparable running times (on CPU) and can take into account periodic boundary conditions as well as distributed point multipoles in a convergent way.
    \item ANKH-FFT is easier to implement than ANKH-FMM, especially in projection for the GPU context while ANKH-FMM, as a IBFMM-based approach, directly benefits from the literature on HPC for this kind of method.
    \item Since ANKH-FFT relies mainly on FFTs, it may be easier to provide portable codes for it (most modern architectures provide efficient FFT implementations).
    \item The drawback, however, is the pre-computation timings of ANKH-FFT that are more costly than ANKH-FMM ones.
\end{itemize}
One of the main motivation behind providing this ANKH approach is to design methods able to efficiently scale on modern HPC architectures and (very) large systems of particles. Hence, it is important to realize that \textit{state-of-the-art} techniques (mainly relying on Local Essential Trees \cite{salmonphd,10.5555/147877.148090,Dubinski_1996,10.1145/2160718.2160740,doi:10.1137/18M1173599}) for parallel computations using FMMs can be used in MD exploiting the ANKH framework. More specifically, Local Essential Trees may also be combined with the ANKH-FFT approach, since local sub-problems can be treated with the latter (only the periodic boundary condition handling has to be slightly adapted). We believe that this last idea, exploiting GPU implementation of ANKH-FFT, can lead to highly scalable MD simulations on modern and hybrid HPC architectures.

\section{Numerical results}
\label{s::numerical}
In this section, we provide numerical experiments aiming at validating our ANKH approach. Follwing Sect. \ref{ss::comparisonfftfmm}, since the two ANKH variants (ANKH-FMM and ANKH-FFT) provide almost the same timing results, the numerical timing results in this section are based on ANKH-FFT only. Our code is sequential and still run on a single Intel(R) Xeon(R) Silver 4116 CPU. We used the  intel oneAPI Math Kernel Library for BLAS calls and FFTW3 for FFT computations. All the code uses double precision arithmetic. The Ewald Summation parameter $\xi$ (see Eq. \ref{eq:ewald}) is always fixed at $0.01$.
 
\subsection{Interpolation errors}
First, we want to verify that our way of treating the point multipole expansions results in admissible errors. To do so, for each partial derivative of $H$ (as in Eq. \ref{eq:Hdef}) in situations corresponding to well-separated pair of cells $(t,s)$, we numerically measured the relative error on the interpolation:
\begin{equation*}
    \partial_\bo{x}^\alpha\partial_\bo{y}^\beta\cc{I}(\bo{x},\bo{y}) := \sum_k\left(\partial_\bo{x}^\alpha S_k[t]\right)\sum_lH(\bo{x}_k,\bo{y}_l)\left(\partial_\bo{y}^\beta S_l[s]\right),
\end{equation*}
\begin{equation}
\label{eq::relerrorforinterpdiff}
    \frac{\displaystyle\mathop{max}_{\substack{\text{ on }100\text{ randomly generated }\\(\bo{x},\bo{y})\in t\times s}}\{\partial_\bo{x}^\alpha\partial_\bo{y}^\beta \left(H(\bo{x},\bo{y}) - \cc{I}(\bo{x},\bo{y})\right)\}}{\displaystyle\mathop{max}_{\substack{\text{ on }100\text{ randomly generated }\\(\bo{x},\bo{y})\in t\times s}}\{|H(\bo{x},\bo{y})|\}},
\end{equation}
where we considered the maximum error, the same set of $100$ test points on both numerator and denominator. The important point is that these errors are given relatively to the maximum magnitude of the (non-differentiated) kernel $H$. This is motivated by the form of our differential operators $\kk{D}_\bo{x}$, $\kk{D}_\bo{y}$ that are linear combinations of partial derivatives, hence influenced by the maximal magnitude of these partial derivative, which is obtained on the kernel evaluation (without derivatives) directly. In addition, because we do not consider scaling by charges, dipoles and quadrupoles moment, and because in our applications the magnitudes of these moments are far smaller than the charges one, the results of Fig. \ref{eq:Hdef} should overestimate the numerical error induced by numerical differentiation on practical cases. Notice that similar results are obtained using more test points (but prohibit the tests for high orders). We tested all possible partial derivatives: since if $|\alpha_0| = |\alpha_1|$ and $|\beta_0| = |\beta_1|$ the errors using $\partial^{\alpha_0}_\bo{x}\partial_{\bo{y}}^{\beta_0}$ and $\partial^{\alpha_1}_\bo{x}\partial_{\bo{y}}^{\beta_1}$ have the same order, we only report one of them in the figure.
\begin{figure}
    \centering
    \includegraphics[width=\linewidth]{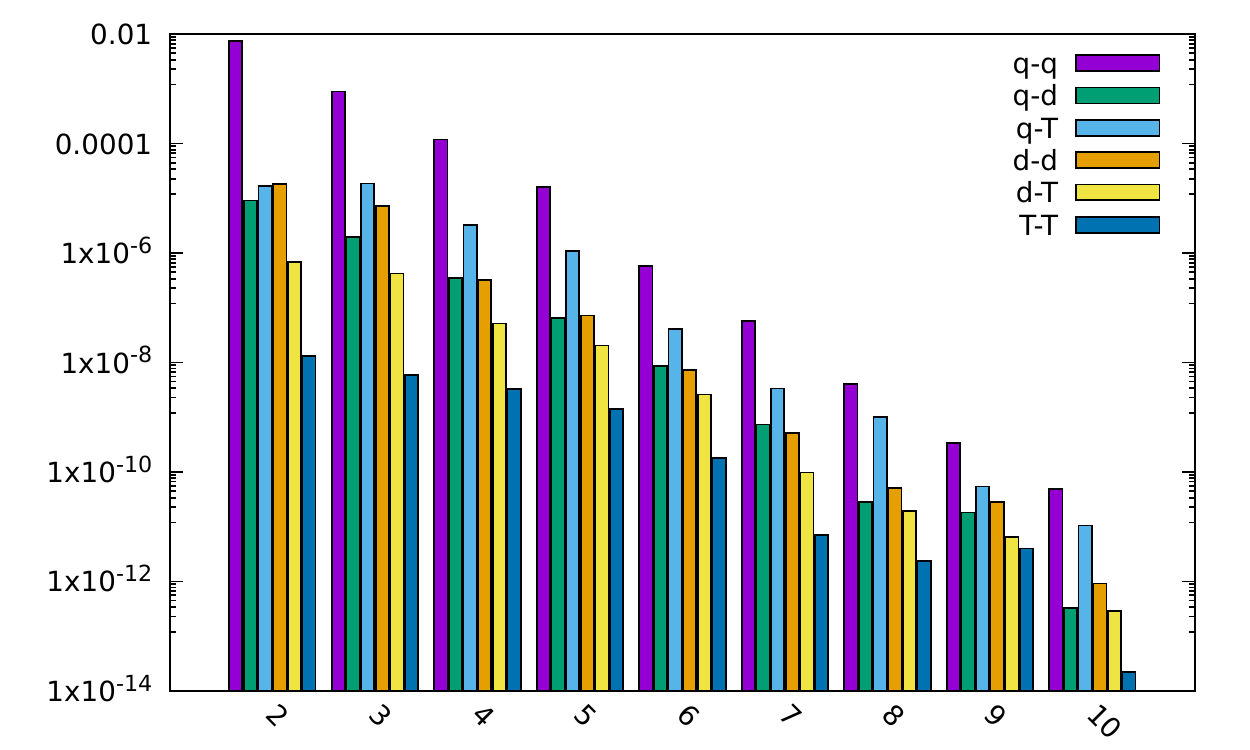}
    \caption{Relative error (according to Eq. \ref{eq::relerrorforinterpdiff}) of kernel interpolation and numerical differentiation. q: charges, d: dipoles, T: quadrupoles. Error types i-j corresponds to derivatives used for i and j applied on kernel $H$.\label{fig:derived_interpolation_error}}
\end{figure}

According to Fig. \ref{eq:Hdef}, the numerical differentiation described in Sect. \ref{ss:numericaldiff} results in errors on derivatives always smaller than the error on the kernel interpolation itself (relatively to this kernel direct evaluation). The numerical error of this interpolation process on the (non-differentiated) kernel follows a $\cc{O}(\nu^L)$ shape, with 1D interpolation order $L$ and $\nu\in (0,1)$), so do all the tested partial derivatives (with the error defined in Eq. \ref{eq::relerrorforinterpdiff}).

As a sum of these ($1+3+6 = 10$, using quadrupole matrix symmetry) partial derivatives, the numerical differentiation error using operators $\kk{D}_\bo{x}$ and $\kk{D}_\bo{y}$ (as in Sect. \ref{s::intro}) should also behaves as $\cc{O}(\nu^L)$. This last point is investigated in Sect. \ref{ss::overallerrorankh}.

\subsection{Overall error in non-periodic configuration}
\label{ss::overallerrorankh}
\begin{figure}
    \centering
    \includegraphics[width=\linewidth]{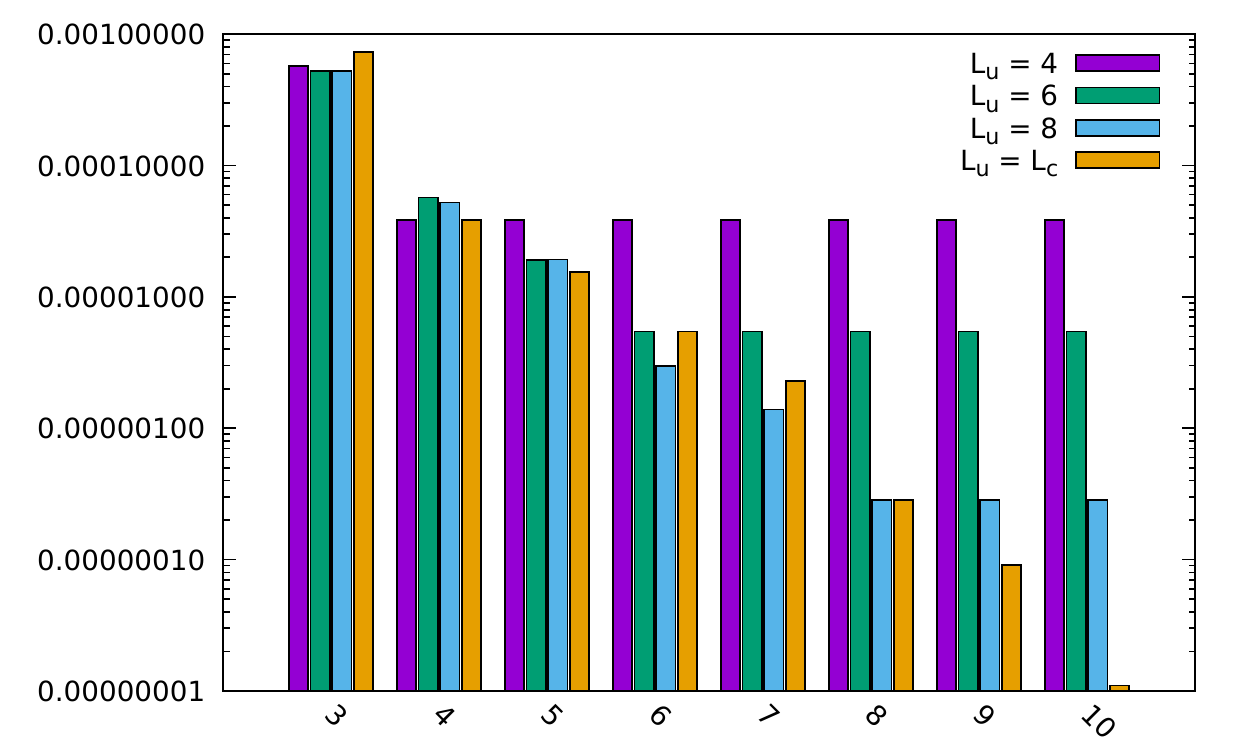}
    \caption{Relative error on water box with 648 atoms (without periodic boundary condition) with respect to the 1D Chebyshev interpolation order $L_c$ and according to various 1D equispaced interpolation order $L_u$.\label{fig:error_no_perio}}
\end{figure}

In this section, we measure the relative error $r_{ANKH}$ of ANKH-FFT results $\cc{E}_{ANKH}$ (relatively to exact solution) on a small water test case without periodic boundary conditions and with respect to 1D (Chebyshev) interpolation order. This relative error is given by
\begin{equation}
\label{eq::relerrnonperio}
    \begin{aligned}
        \cc{E}_{exa} &= \frac{1}{2}\sum_{\bo{x}\in X}\sum_{\bo{y}\in Y}\kk{D}_{\bo{x}}\kk{D}_{\bo{y}}\left(\frac{1}{|\bo{x}-\bo{y}+\bo{t}|}\right),\\
        r_{ANKH} &= \frac{|\cc{E}_{exa} - \cc{E}_{ANKH}|}{|\cc{E}_{exa}|}.
    \end{aligned}
\end{equation}
According to Rem. \ref{rem::compression}, the 1D interpolation order $L_u$ on equispaced grids does not have to fit with the Chebyshev interpolation one $L_c$, hence we considered different $L_u$ for the tested $L_c$ in order to check the compression error that may be induced by choosing $L_u < L_c$. Results are provided in Fig. \ref{fig:error_no_perio}.

First, we observe that the error geometrically decreases with respect to the 1D Chebyshev interpolation order, provided that the equispaced one equals this order. Hence, the overall ANKH-FFT approach numerically converges in accordance with the theory. According to the geometric convergence, when fixing the equispaced order to $4$, the relative error does not decrease under $\approx 10^{-5}$, while it stagnates at $\approx 10^{-6}$ and $\approx 10^{-7}$ for equispaced orders fixed to $6$ and $8$ respectively. Moreover, when fixing to a constant the equispaced 1D interpolation order, the convergence is still observed up to the precision allowed by the latter (and then stagnates). This thus validates the trick we presented in Rem. \ref{rem::compression}: we can further compress the Chebyshev interpolation process using possibly smaller equispaced interpolation orders.

\subsection{Overall error in periodic boundary configuration}
\label{ss::overallperioankh}
\begin{figure}
    \centering
    \includegraphics[width=\linewidth]{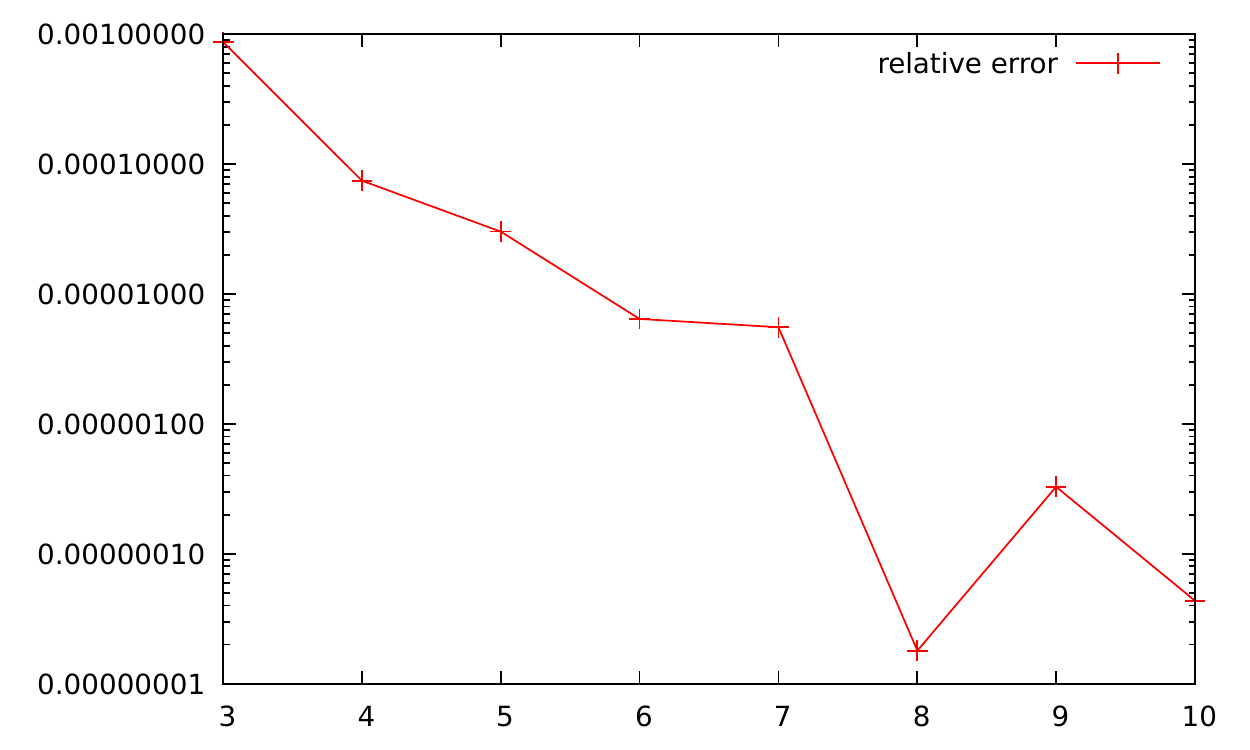}
    \caption{Relative error of ANKH-FFT de with respect to the 1D Chebyshev interpolation order.\label{fig:ovperioankh}}
\end{figure}

The tested relative error is the same than in Eq. \ref{eq::relerrnonperio} except that periodic boundary conditions are used to $\cc{E}_{exa}$, fitting with the target quantity defined in Eq. \ref{eq_qmutheta_coulomb}:
\begin{equation}
    \begin{aligned}
        \cc{E}_{exa} &= \frac{1}{2}\sum_{\bo{t}\in 2r_B\bb{Z}^3}\sum_{\bo{x}\in X}\sum_{\bo{y}\in Y}\kk{D}_{\bo{x}}\kk{D}_{\bo{y}}\left(\frac{1}{|\bo{x}-\bo{y}+\bo{t}|}\right),\\
        r_{ANKH} &= \frac{|\cc{E}_{exa} - \cc{E}_{ANKH}|}{|\cc{E}_{exa}|}.
    \end{aligned}
\end{equation}
The corresponding results are provided in Fig. \ref{fig:ovperioankh}, in which we considered $31^3-1$ images of the simulation box, an 1D equispaced interpolation order equal to the Chebyshev one $L_c$ and interpolation order for the far images handling (see Sect. \ref{ss::effperiofft}) equal to $L_c+1$ (which was empirically fixed and may not be optimal). The test case is the same than in Sect. \ref{ss::overallerrorankh}. This choice of image number is also empirical: we observed in our tests that such parameter allowed to fully converge in most cases.

According to these errors, our ANKH approach converges in the periodic boundary condition configuration, still geometrically (hence as the non-periodic case of Sect. \ref{ss::overallerrorankh}). The relative error value for $L_c=8$ is surprisingly low, but we believe that this is just a particular case on which more digits are luckily correct than expected in this particular test case and using this particular parameter choice. This validates our ANKH approach when using periodic boundary conditions.

Unfortunately, due to the cost of direct computation of exact reference solution using large number of far images, we cannot directly check our code with exact solution when considering important number of atoms. However, to emphasize the conclusion in this section, as discussed (and done) in Sect. \ref{ss::performancetest}, we can verify that our ANKH approach is able to reach the same accuracy than provably convergent code on larger cases. Notice that we also retrieve the geometrical convergence under periodic boundary conditions when using few number of images (to perform the direct computation) on test cases with $\approx 100000$ particles (hence, since the conclusions are the same, we do not report them in this paper).

\subsection{Performance test}
\label{ss::performancetest}
In this section, we compare the performance of both ANKH-FFT and Smooth Particle Mesh Ewald (SPME) using Tinker-HP on various systems. Their size range from a small water box of 216 water molecules (648 atoms) up to the Satellite Tobacco Mosaic Virus in water (1066624 atoms), covering typical scale of systems of interest in biochemistry even if both methods could naturally handle larger cases. In more details these are:
\begin{itemize}
    \item water box of 648 atoms,
    size 18.643$^3$ $\mathring{A}$,
    \item dhfr protein in water, 23558 atoms, size 62.23$^3$ $\mathring{A}$,
    \item water box of 96000 atoms, size 98.65$^3$ $\mathring{A}$,
    \item water box of 288000 atoms, size 142.17$^3$ $\mathring{A}$,
    \item water box of 864000 atoms, size 205.19$^3$ $\mathring{A}$,
    \item Satellite Tobacco Mosaic Virus in water, 1066624 atoms, 223.0$^3$ $\mathring{A}$.
\end{itemize}
Most of these test cases can be found in the Tinker-HP git repository: \url{https://github.com/TinkerTools/tinker-hp/tree/master/v1.2/example}

Let us recall that we compare sequential execution of both applications without including precomputations (see Sect. \ref{ss::comparisonfftfmm}). We use Tinker-HP results with default SPME parameters as a reference, which we observed to be of order $10^{-4}$ compared to pure Ewald summation. Thus, we tuned the ANKH-FFT parameters to reach similar relative errors. We considered, for these tests, a 1D equispaced interpolation order equal to $2$. 
\begin{figure}
    \centering
    \includegraphics[width=\linewidth]{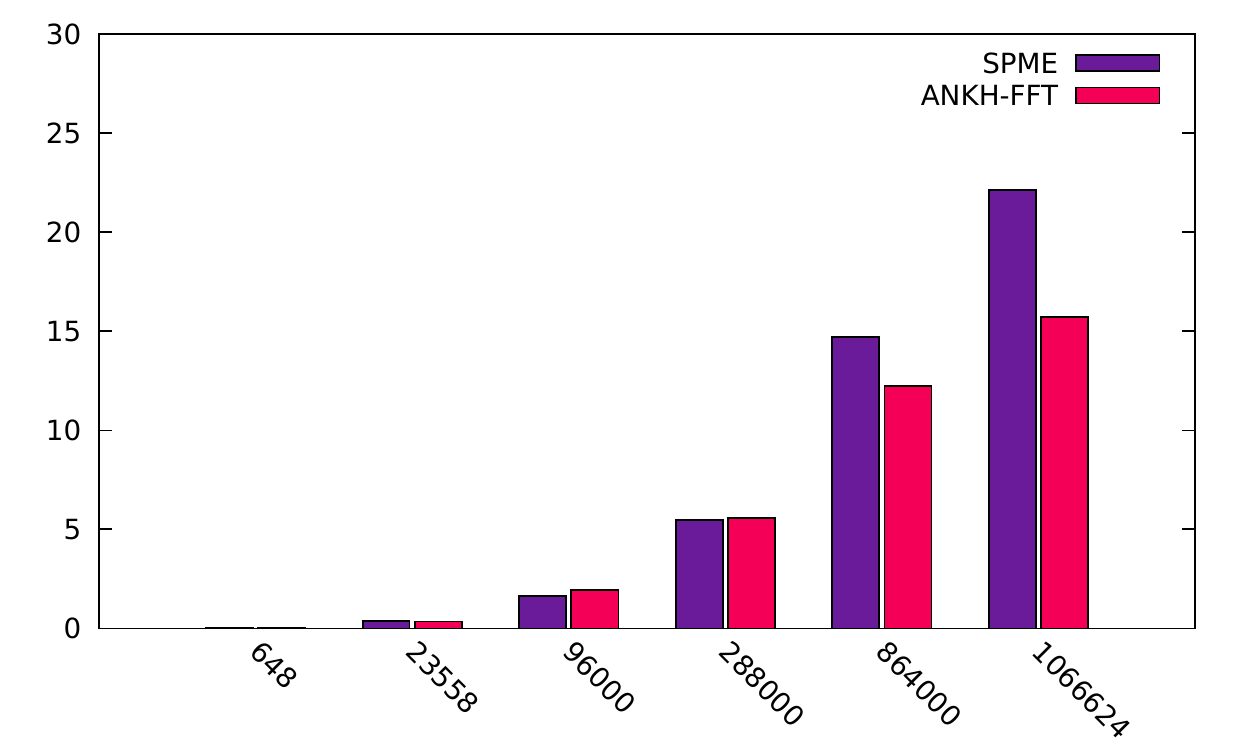}
    \caption{Comparison of application timings (in second) for ANKH-FFT (red) and SPME (purple) with regard to the number of atoms in $B$.\label{fig:ANKHvsTINKER}}
\end{figure}

On small test cases (between few hundreds and few thousands), ANKH-FFT and SPME show similar timings. Medium test cases appears to be the most favorable to SPME here, where it perform $20\%$ faster than ANKH-FFT (for 96000 atoms). Notice that this number of particle is quite small for hierarchical methods and that performance still are very similar (same magnitude order each time), which allows us to claim that our method almost perform the same way than \textit{state-of-the-art} approaches on this type of particle distributions. On larger test cases, ANKH-FFT outperforms SPME, providing timings up to $40\%$ smaller than the last.

Hence, ANKH-FTT shows important performance on wide range of particle distribution size and for various types of typical systems arising in biochemistry. This direct code performance comparison however does not allow to conclude on the interest of the optimizations we provided. We thus focus on ANKH-FFT on Sect. \ref{ss::splittimings}.

\subsection{Split timings and details}
\label{ss::splittimings}
\begin{figure}
    \centering
    \includegraphics[width=\linewidth]{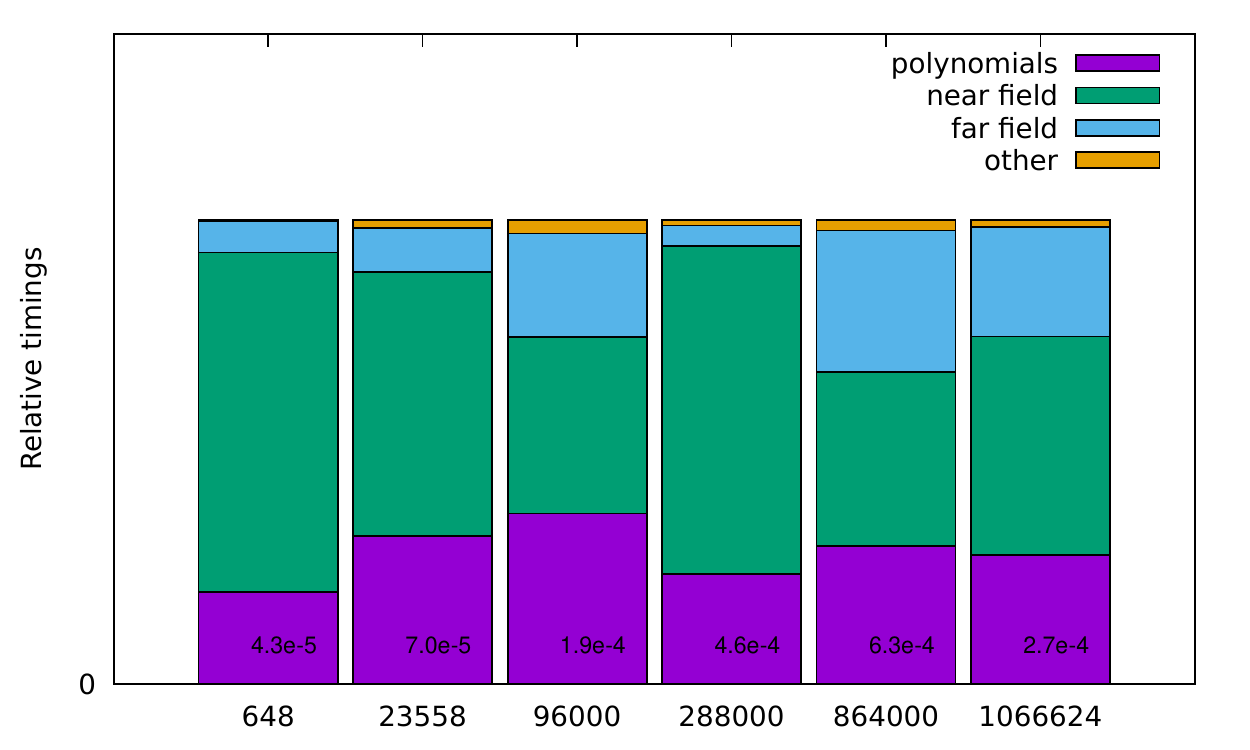}
    \caption{Relative timings (with respect to the entire running time of each test case) of the different steps of ANKH-FFT application. Relative error (compared to the results of Tinker-HP) are indicated in black and the number of atom in each corresponding system corresponds to the abscissa.\label{fig:rel_timings_overall}}
\end{figure}

In this section we provide additional details on the tests of Sect. \ref{ss::performancetest}. First, we present in Fig. \ref{fig:rel_timings_overall} split timings in order to compare the running times of the various evaluation steps. There mainly are three of them:
\begin{itemize}
    \item the computation of interpolation polynomials at leaves and transformation to modified charges,
    \item the far field computation (using FFT),
    \item the near field computation (using direct evaluation).
\end{itemize}

Clearly, the most time consuming part is dedicated to near field computations. This step, as well as the polynomial handling, has the strong advantage of being a local one. In addition, tree structure localizes the particles in an efficient way, avoiding the computation of neighbor lists. We thus expect these two step to fully benefit from parallelization.

Since the far-field computation time does not exceed $\approx 6\%$ of the overall application time, these results tend to validate our fast handling of far interactions. Notice that this far field still involves most of the interactions (at least for medium and large test cases). However, these relatively small timings are mainly due to the important recompression, using 1D equispaced interpolation order equal to $2$. This was sufficient to reach the Tinker-HP SPME implementation precision (see relative errors on Fig. \ref{fig:rel_timings_overall}). For higher required precisions, the balance between near and far fields is better preserved in overall timings. Nevertheless, this recompression allows to strongly mitigate the hierarchical methods prefactor.

\begin{figure}
    \centering
    \includegraphics[width=\linewidth]{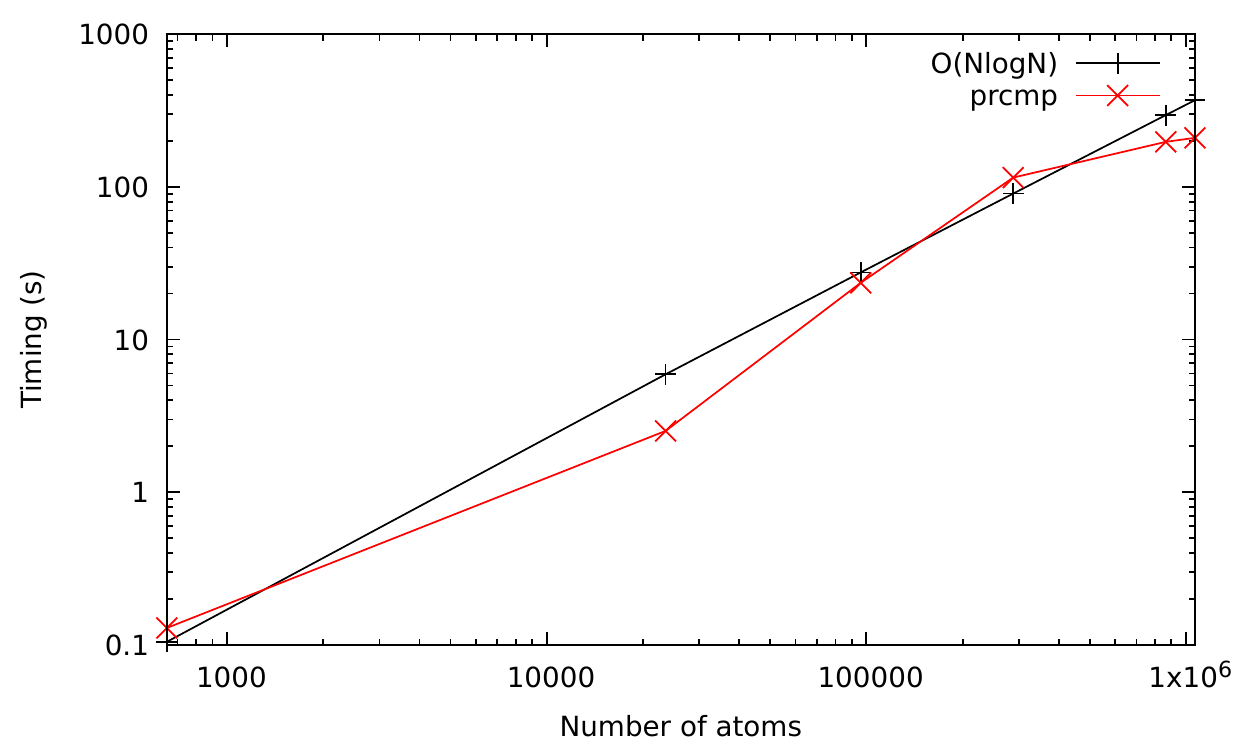}
    \caption{Precomputation timings with respect to number of atoms for fixed $10^{-4}$ output precision and using $31^3-1$ images of the simulation box $B$ (in red), compared to $\cc{O}(NlogN)$ curve (in black)\label{fig:precomp_comp}}
\end{figure}

As stated in Sect. \ref{ss::ankhfft_complexity}, ANKH-FFT precomputation complexity is linearithmic. These costs are negligible in applications according to Sect. \ref{ss::comparisonfftfmm}. We still present in Fig. \ref{fig:precomp_comp} the precomputation timings of ANKH-FFT on the different test cases of Sect. \ref{ss::performancetest}.

According to these timings, the precomputation step behaves as our theoretical linearithmic estimates. This is a consequence of the interpolation used for the handling of far field images since we also compared to a non-interpolated version that was already not able to run in medium test cases. In terms of magnitude order, ANKH-FFT precomputation timings are up to few hundred seconds on larger test cases. Regarding the entire cost of MD simulation, since these precomputations are done only once, the conclusion proposed in Sect. \ref{ss::comparisonfftfmm} is numerically verified.

\section{Conclusion}
In this paper, we introduced ANKH, a new polynomial interpolation-based accelerated Ewald Summation for energy computations in molecular dynamics. Our method mainly exploits the (absolutely convergent) real part of the Ewald Summation in order to express the energy computation problem as a generalized $N$-body problem on which fast kernel-independent hierarchical method can perform. We then proposed two fast approaches, the first one based on IBFMM (accordingly named ANKH-FMM) will be well-suited for distributed memory parallelization and the second, built on a fast matrix diagonalization using FFTs (named ANKH-FFT) will be a good candidate for GPU computations. ANKH-FMM has a linear complexity while ANKH-FFT has a linearithmic one but practically behaves as a linear complexity approach. We demonstrated the numerical convergence of our interpolation-based method and we compared our implementation to a Smooth Particle Mesh Ewald production code (Tinker-HP), exhibiting comparable performance on small and medium test systems and superior performances on large test cases. A public code is provided at \url{https://github.com/IChollet/ankh-fft} whereas a version of it will also be available within a future Tinker-HP code release (\url{https://github.com/TinkerTools/tinker-hp}).

We believe that ANKH is an interesting approach for MD on modern HPC architectures, especially because it theoretically solved the scalability issues of Particle-Mesh methods in distributed memory context. However, our current implementation does not use shared or distributed memory parallelism and our intuition still has to be verified in practice. Hence, our priority is now to develop a scalable ANKH-based code (MPI + GPU) and to test it on early exascale architectures.

Among other short-term direct perspectives of the work we presented here, we are considering the (almost) straightforward extension to forces computations. Indeed, the main difference between the algorithm we presented and a general force algorithm belong in the mutual interactions that cannot be exploited the same way. We would also like to look at application of ANKH-FFT to other scientific fields on which $N$-body problems on quasi-uniform particle distributions have to be solved. 

Considering long term perspectives, the adaptation of the algorithms we presented to non-cubic simulation boxes may also find applications in material science.

\section*{Acknowledgements}
This work has also received funding from the European Research Council (ERC) under the European Union's Horizon 2020 research and innovation program (grant agreement No 810367), project EMC2 (JPP).
\bibliographystyle{plain}
\bibliography{bibli}

\end{document}